\documentclass{article}

\usepackage[pdftex,a4paper, height=230mm,width=164mm]{geometry} 
\usepackage[margin=15pt,font=small,labelsep=period,labelfont=sc]{caption} 
\usepackage[bf,compact,small,pagestyles]{titlesec}

\usepackage{stix2}

\usepackage{mathrsfs} 
\usepackage{microtype}

\usepackage[absolute,overlay]{textpos} 

\usepackage{graphicx}
\usepackage[title]{appendix}

\usepackage{mathtools}
\usepackage{amsthm}
\usepackage{stmaryrd} 		

\let\mathbb\undefined 
\usepackage{BOONDOX-ds} 

\numberwithin{equation}{section}  

\newtheorem{theorem}{\textbf{Theorem}}[section]
\theoremstyle{remark}
\newtheorem{remark}[theorem]{\textit{Remark}}
\newtheorem{assumptions}[theorem]{\textit{Assumptions}}

\newcommand{\e}{\textrm e}
\renewcommand{\i}{\textrm i}
\newcommand{\Real}{\operatorname{Re}}
\newcommand{\RR}{\mathbb R}
\let\bs\boldsymbol

\newcommand{\DeltaT}{\Delta_\text{T}}
\newcommand{\nablaT}{\nabla_\text{T}}

\newcommand{\holdall}{D}
\newcommand{\dif}[1]{\textit d#1}

\newcommand{\dc}[1]{\mbox{\fontdimen2\textfont1=1pt$\partial#1\null$\fontdimen2\textfont1=0pt}}

\newcommand{\der}[3]{\ifcat#31{
                \ifnum#3=1 {\frac{\textit d#1}{\textit d#2}}
                \else{\frac{d^#3#1}{\textit d#2^#3}}
                \fi}
        \else{\frac{d^#3#1}{\textit d#2^#3}}
        \fi}
\newcommand{\pder}[3]{\ifcat#31{
                \ifnum#3=1 {\frac{\dc{#1}}{\dc{#2}}}
                \else{\frac{\partial^#3#1}{\dc{#2}^#3}}
                \fi}
        \else{\frac{\partial^#3#1}{\dc{#2}^#3}}
        \fi}

\title{A better compression driver?\\
CutFEM 3D shape optimization  taking viscothermal losses\\  into account}
\author{Martin Berggren\textsuperscript1, Anders Bernland\textsuperscript1, André Massing\textsuperscript2, \\
Daniel Noreland\textsuperscript3, and Eddie Wadbro\textsuperscript4\\[2pt]
\small\itshape
\textsuperscript1Department of Computing Science, Umeå University, Sweden\\[-2pt]
\small\itshape
\textsuperscript2Department of Mathematical Sciences, Norwegian University of Science and Technology, Norway\\[-2pt]
\small\itshape
\textsuperscript3The Forestry Research Institute of Sweden, Uppsala, Sweden\\[-2pt]
\small\itshape
\textsuperscript4Department of Mathematics and Computer Science, Karlstad University, Sweden
}
\date{\normalsize March 12, 2024}

\begin{document}

\maketitle

\begin{abstract}

The compression driver, the standard sound source for midrange acoustic horns, contains a  cylindrical compression chamber connected to the horn throat through a system of channels known as a \textit{phase plug}.
The main challenge in the design of the phase plug is to avoid resonance and interference phenomena. 
The complexity of these phenomena makes it difficult to carry out this design task manually, particularly when the phase-plug channels are radially oriented.
Therefore, we employ an algorithmic technique that combines numerical solutions of the governing equations with a gradient-based optimization algorithm that can deform the walls of the phase plug.
A particular modeling challenge here is that viscothermal losses cannot be ignored, due to narrow chambers and slits in the device. 
Fortunately, a recently developed, accurate, but computationally inexpensive boundary-layer model is applicable. 
We use this model, a level-set geometry description, and the Cut Finite Element technique to avoid mesh changes when the geometry is modified by the optimization algorithm. 
Moreover, the shape calculus needed to compute derivatives for the optimization algorithm is carried out in the fully discrete case.
Applying these techniques, the algorithm was able to successfully design the shape of a set of radially-directed phase plugs so that the final frequency response surprisingly closely matches the ideal response, derived by a lumped circuit model where wave interference effects are not accounted for.
This result may serve to resuscitate the radial phase plug design, rarely used in today’s commercial compression drivers.

\end{abstract}
\begin{textblock*}{\textwidth}(1.5cm,1.5cm) 
\textit{Preprint of submitted manuscript}
\end{textblock*}


\bigskip

\section{Introduction}

We present  an engineering design case study  that aims to demonstrate  how recent progress, both in mathematical modeling and design optimization techniques, may potentially revive a  design concept rarely used  due to difficulties in providing simple explicit design guidelines.
The mathematical modeling progress concerns acoustic viscothermal (also called thermoviscous)  losses, and the improved design optimization technique involves a novel shape calculus approach for CutFEM discretizations  in combination with level-set geometry descriptions. 
The target application we consider is the \textit{compression driver}, a standard sound source for the acoustic horns that are used in public address systems aimed at large auditoria and outdoors.
In particular, we will optimize the device in the rarely-used \textit{radial} orientation of the so-called phase plugs.
We will now shortly introduce these improved techniques, introduce the device under consideration, and outline the content of the rest of the contribution.

\subsection{Viscothermal design optimization}

The mathematical modeling of viscothermal losses is of importance for acoustic devices in which the acoustic waves interact with large areas of solid surfaces, for instance in waveguides, acoustic liners, and in geometrically narrow devices.
One example of a narrow device is the compression driver considered here;  other examples are miniaturized devices such as headphones, microphones, mobile phones, and hearing aids.

It is only recently that numerical optimization of acoustic devices has been carried out using methods that can account for viscothermal losses. 
To a high computational cost, particularly for the three-dimensional cases we consider here, such losses can within the linear regime be accurately modeled using the linearized, compressible Navier--Stokes equations.
In addition to the cost, a  problematic issue of using these equations in the the context of design optimization is the extreme scale separation between the viscothermal and wave propagation effects. 
The viscothermal losses are typically concentrated in exceedingly thin boundary layers close to solid walls.
For instance, for air in the audio regime 20\,Hz--20\,kHz, the thickness of these boundary layers is smaller than the wavelength by a factor of $10^{-5}$--$10^{-3}$~\cite{BeBeNo18}.
Thus, to achieve reasonable computing times, it is necessary to employ aggressive mesh refinement strategies along solid boundaries, preferably in the form of highly anisotropic meshes: thin in the wall-normal direction and stretched-out in the tangential direction.
The management of such meshes in a shape-optimization context is challenging.
 
However, the scale separation between wave and viscothermal effects can also be exploited to arrive in a simplified model of lower computational cost.
Such a model, the sequential linear Navier--Stokes model, was developed by Kapinga~\cite{Ka10}.
As for the full Navier--Stokes equations, this simplified model unfortunately still requires accurate spatial resolution of the boundary layers.
Noguchi \& Yamada~\cite{NoYa22} use the Kapinga model to optimize a sound-absorbing pipe in 2D axial symmetry, employing isotropic mesh refinements based on approximate computations of the signed distance function.
The design sensitivity analysis utilized the adjoint-equation approach; however, the contributions from the viscothermal terms in the final gradient expressions were ignored in the optimization algorithm.

One way to circumvent the need for resolution of the boundary layers is to use a boundary-element formulation, which does not require a mesh inside the domain. 
Such an approach was used by Andersen, Henriquez, and Aage~\cite{AnHeAa19,AnHeAa23} to optimize the shape of quarter-wave and Helmholtz resonators in 2D axial symmetry.
Their viscothermal model is similar to the Kapinga model, but the resulting Helmholtz and diffusion equations are solved using boundary-element techniques.
Although the use of the boundary-element method reduces the spatial dimensionality by 1 and circumvents the need for volume meshes, the computational cost for solving the system of equations appears to be somewhat high.

The low reduced frequency model~\cite{Ti75} is an even simpler and faster viscothermal model that can be used for special geometries, such as narrow slits and tubes.
Christensen~\cite{Cr17}  used such a model for material-distribution topology optimization, where the cross-section shape of a narrow pipe was designed to maximize the viscous losses. 

A completely different class of viscothermal models is based on boundary-layer theory, such as the model devised by three of the authors to the current contribution~\cite{BeBeNo18}.
This model uses the standard, isentropic Helmholtz equation in the interior of the domain and models the viscothermal losses through a generalized impedance (a so-called Wentzell) boundary condition on solid boundaries.
The computational cost is essentially the same as for lossless acoustics, and the accuracy is comparable with the full linearized Navier--Stokes equation as long as boundary layers from opposing walls do not overlap and as long as the wall curvatures are not too extreme.
In fact, the model seems to work surprisingly well even at really small scales, such as in the modeling of porous materials~\cite{CoDaMaBaBe20}.
This viscothermal model has been used for design optimization only a few times so far and only in 2D.
The model was used by Tissot, Billard, and Garbard~\cite{TiBiGa20} to optimize for maximum absorption the interior shape of a Helmholtz resonator, thought of as an element to be periodically repeated to form an acoustic liner.
Similarly to what we do here, the geometry was specified using a level-set function, and a CutFEM approach was used for discretization, whereas, in contrast to our approach, a fully discrete shape calculus was not attempted in their contribution.
Dilgen, Aage, and Jensen~\cite{DiAaJe22} considered topology optimization of a hearing  instrument, utilizing a vibroacoustic model for which the viscothermal losses were modeled using the boundary-condition approach.
They also employed a level-set-function geometry description in combination with a CutFEM discretization approach.
Parts of the calculations involved in the sensitivity analysis were carried  out approximately using finite differences. 
Recently, Mousavi et al.~\cite{MoBeWa23,MoBeHaWa24} introduced the boundary-layer model in a completely different optimization scenario, using so-called material-distribution (or density-based) topology optimization.
In this scenario, the geometry is represented as a varying density-like coefficient in the governing equations.
Such methods can produce designs with arbitrary topological complexity to the price of producing stair-cased boundaries and needing very fine meshes to accurately approximate curved boundaries.

\subsection{The target application}

The most common sound source used in loudspeakers is the moving-coil dynamic transducer.
When the diaphragm of such a driver is oscillating in free air, only a small fraction of its vibrational energy is converted to propagating sound waves due to the mismatch of acoustic impedance between air and the solid material in the diaphragm. 
Already in the 1920s, Hanna \& Slepian~\cite{HaSl77} realized that better impedance matching could be achieved by placing the diaphragm in a  cylindrical compression chamber with a narrow outlet attached to the throat of an acoustic horn.
It was soon realized that the way in which the outlet from the compression driver is connected to the horn throat is important for acoustic performance; this part of the device is now known as the \textit{phase plug}. 
The electrodynamic motor, the compression chamber, and the phase plug are typically mounted together to form a \textit{compression driver},  which to this day is universally used  to feed acoustic horns at mid to high frequencies in public address systems.
The book by Kolbrek and Dunker~\cite{KoDu19} contains a recent and extremely comprehensive discussion on the history as well as theory of acoustic horns and drivers.

A main challenge in the design of the compression driver is to avoid resonance and interference phenomena, as they cause an uneven frequency response and may aggravate nonlinearities. 
A first issue is radial modes appearing in the compression chamber.
In most commercially available devices, the compression chamber contains several outlets oriented in the annular direction, relying on the ingenious design rules first introduced by Smith~\cite{Sm53}. 
Using  such  \textit{Smith slits}, that is, placing $N$ circumferential channels following the guidelines, $N$ radial modes in the compression chamber can be suppressed.
Dodd and Oclee-Brown~\cite{Do03,OB12} modified Smith’s guidelines to curved and Voishvillo~\cite{Vo16} to ring-shaped diaphragms.

However, the phase plug that connects the slits to the horn throat can create additional resonances through wave interference between the channels.
A reasonable guideline for the design of the channels in the phase plug is that these should be of equal path length.
However, this requirement still does not prevent internal resonances from occurring.
Bezzola~\cite{Be18} briefly reported results of optimizing the channels using the commercial Comsol Multiphysics software (likely using the lossless Helmholtz equation).
In a previous contribution, three of the authors of the current article applied numerical shape optimization to the design of the channels in the phase plug of a compression driver equipped with Smith slits~\cite{BeWaBe19}.
We concluded that the shape optimized channels indeed improved the frequency response compared to the initial straight-channel design, but that multiple resonances unfortunately remained also in the optimized device.
Another conclusion from the study was that viscothermal boundary-layer losses are quite prominent, due to the wave propagation being confined in a very narrow geometry.
To lower the complexity and computational cost of the calculations, we carried out the optimization neglecting these losses but evaluated the performance on the final geometry accounting for the viscothermal losses.
However, thanks to recent modeling advances~\cite{BeBeNo18}, we will here be able to include the modeling of the losses also in the optimization step.

Here we will consider the case when the outlets from the compression chamber are oriented in the \textit{radial} direction.
This design choice is rare in commercial units, and only little analysis has been carried out.
The most complete discussion is likely the chapter on radial outlets in the thesis of Oclee-Brown~\cite{OB12}.
In order to fit to typical geometries found in compression drivers, Oclee-Brown studies the wave propagation in the compression chamber by expanding the solution to the Helmholtz equation for the acoustic pressure in spherical coordinates.
The analysis suggests that the number of outlets should be large enough to suppress circumferential modes in the chamber, and that radial modes can be suppressed by a slight modification of constant annular-width channel outlets.

Concerning the design of a phase plug that connects  radially-oriented outlets into a waveguide to be fed into a horn, there seem to be no  guidelines at all based on acoustic analysis in the literature.
Here we will utilize gradient-based shape optimization for this task, based on a further development of the method we used for designing phase plugs for circumferential Smith slits~\cite{BeWaBe19}.
The method relies on a level-set description of the geometry, which allows complex shapes to be treated in the optimization.
In previous projects, the use of a level-set description has allowed unexpected and nontrivial shapes to appear  in acoustic shape optimization studies~\cite{BeWaBe17,BeWaBe19}.
We combine the level-set geometry description with a fictitious-domain method known as CutFEM~\cite{BuClHaLaMa15} or XFEM.
That is, the boundary of the computational domain, here the phase-plug wall, is allowed to pass arbitrarily through the computational mesh, which means that only partial contributions to the governing equations come from the finite-element mesh elements that are intersected by the boundary.
Moreover, boundary conditions for the walls of the phase plug have to be assigned to the boundary surfaces that cut through the elements.
As opposed to the previous study~\cite{BeWaBe19}, the calculations need here to be carried out in three space dimensions. 
Previously, CutFEM/XFEM approaches have successfully been used for the shape optimization of, for instance, elastic structures~\cite{DuysinxMiegroetJacobsEtAl2006,BurmanElfversonHansboEtAl2018,BretinChapelatOuttierEtAl2022},
vibroacoustic problems~\cite{DiAa21,DiJeAa21,DiAa24}, and flow patterns~\cite{VillanuevaMaute2017,DokkenJohanssonMassingEtAl2020}.


\subsection{Overview of article}

In \S\,2 and appendix\,\ref{a:lumpedmodel}, we review the acoustic mechanism of the compression driver and derive an expression that gives an upper limit of the performance, an expression that will be used in our most successful optimization formulation.
Modeling details, the optimization formulation, and details of the discretization for our target application are introduced in \S\,4--6. 
A major difference with the previous study~\cite{BeWaBe19} is that we here will consider viscothermal losses also during the optimization, based on the model discussed above~\cite{BeBeNo18}. 
The viscothermal boundary condition needs to be supplied to the walls of the phase plug subject to design, which significantly complicates the shape calculus.
For the first time, to the best of our knowledge, this calculus will be carried out in the fully discrete case, taking into account that the derivatives of the acoustic pressure as well as the normal field of the boundary are in general discontinuous between elements.
For this calculus, we will rely on formulas proven in a recent publication~\cite{Be23}.
The calculus is detailed in appendix\,\ref{a:shapecalc} and the formulas are summarized in \S\,7.
We present the results of the phase plug optimization in \S\,8, results which seem to indicate favorable acoustic properties of properly shaped radially-oriented phase plugs, as discussed in the concluding \S\,9. 

\section{The compression driver mechanism}\label{Slumpedmodel}

To recognize the issues faced when designing a compression drive phase plug, it is necessary first to appreciate the fundamental  mechanism of the compression driver.
In any loudspeaker transducer, a vibrating diaphragm induces the  acoustic waves that will be transmitted to the surroundings.
The basic function of  the compression driver is to generate a high sound pressure level while only requiring a small displacement of the diaphragm.
Hanna \& Slepian~\cite{HaSl77} are generally recognized to be the first in 1924 to describe and analyze such a device.

Under a number of strong assumptions, it is possible to devise a lumped-parameter model for the acoustic pressure generated by the movement of a diaphragm in a model compression driver.
We will see that only two geometric properties of the device, its compression ratio and the depth of the compression chamber, govern its performance under these assumptions.
These two quantities constitute the fundamental parameters to be chosen by the designer.

\begin{figure}[b]\centering
\includegraphics{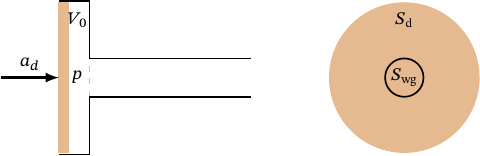}
\caption{A conceptual cylindrical compression driver. 
		Cross-section view  through axis (left) and orthogonal to axis~(right).}\label{f:chamber}
\end{figure}

Consider the  simplified model of a compression driver illustrated in figure\,\ref{f:chamber}.
A diaphragm of cross sectional area $S_\text{d}$ is positioned in a cylindrical chamber with rest volume $V_0$ and with an outlet attached to a waveguide with cross sectional area $S_\text{wg}$.
We make the following idealized assumptions.
\begin{assumptions}\label{ass}
\begin{enumerate}
\item[] 

\item[(i)]\label{asspiston}
The diaphragm moves axially as a solid piston given a time-harmonic movement with acceleration amplitude~$a_d$.

\item[(ii)]\label{assPconst}
The pressure is the same  at each point within the compression chamber 

\item[(iii)]\label{assisentrop}
The compression within the chamber is isentropic, and the air is an ideal gas.

\item[(iv)]\label{asslinear} 
Linearity of the acoustic quantities.

\item[(v)] The waveguide is perfectly terminated; that is, the waves traveling to the right in the waveguide will be perfectly absorbed.

\end{enumerate}
\end{assumptions}

Due to the time-harmonic assumption, the acoustic pressure in the compression chamber and the acceleration of the diaphragm, as a function of time, can be written
\begin{equation}
P(t) = \Real\bigl(p\, \e^{\i\omega t}\bigr), \qquad
A_d(t) =  \Real\bigl(a_d\,\e^{\i\omega t}\bigr),
\end{equation}
where $\omega$ is the angular frequency of the movement and $p$ and $a_d$ are  complex amplitudes.
Under assumptions~\ref{ass}, it can be shown that the acoustic pressure and diaphragm acceleration amplitudes are related through the equation
\begin{equation}\label{e:lumpedformula}
k\left(-dk + \frac{\i}\kappa\right ) p = \rho_0 a_c.
\end{equation}
The geometric parameters of formula~\eqref{e:lumpedformula} are
\begin{equation}
d = \frac{V_0}{S_\text d}, \qquad \kappa = \frac{S_\text d}{S_\text{wg}},
\end{equation}
that is, the compression chamber depth $d$ in the simple geometry of figure\,\ref{f:chamber} and the compression ratio $\kappa$.
Moreover,  $\rho_0$ is the static density of air and $k = \omega/c_0$ the wave number, in which $c_0$ is the speed of sound. 

Formula~\eqref{e:lumpedformula} is likely no surprise for the acoustician and could be extracted, for instance, from the much more comprehensive lumped-parameter performance analysis found in Kolbrek \& Dunker's  monograph~\cite[\S\,18.4.4]{KoDu19}.
However,  the precise formulation given above will be useful in the following, and we are not aware of any  source containing this exact formula, which is why we provide a short derivation in~appendix\,\ref{a:lumpedmodel}.

\begin{figure}
\centering
\includegraphics{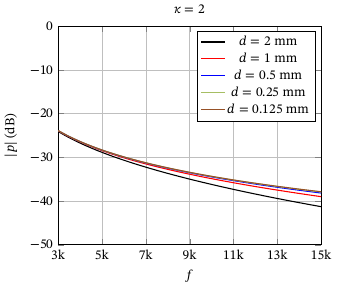}
\qquad\qquad
\includegraphics{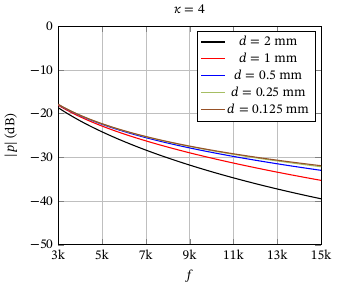}
\\
\includegraphics{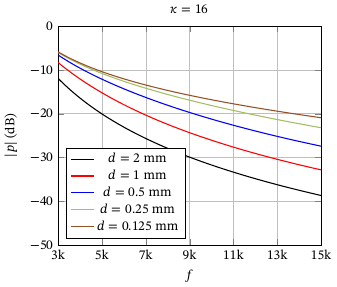}
\qquad\qquad
\includegraphics{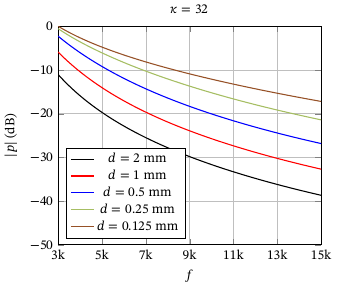}
\caption{The frequency response, given a fixed diaphragm acceleration, of the lumped model for various compression ratios $\kappa$ and chamber depths $d$. 
(Plots scaled by the maximum sound pressure of all graphs.)}\label{f:lumpedfresponse}
\end{figure}

To appreciate the effects of the compression chamber on the frequency response, we will consider the pressure amplitude, which by expression~\eqref{e:lumpedformula} will be
\begin{equation}\label{e:abslumpedformula}
\lvert p\rvert = \frac{\rho_0\lvert a_\text{c}\rvert}{k(d^2k^2 + 1/\kappa^2)^{1/2}}.
\end{equation}
Figure\,\ref{f:lumpedfresponse}  shows the frequency response of the pressure amplitude, according to formula~\eqref{e:abslumpedformula} with constant $\lvert a_\text{c}\rvert$, for a variety of compression ratios and chamber depths. 
We note that   an increased compression  ratio has a dramatically positive effect on the delivered sound pressure level, and that an increase in compression ratio should be accompanied by a narrowing  chamber depth in order not to worsen the dropoff for higher frequencies.

\begin{remark}\label{r:BiskopBrask1}
Note that no real compression driver unit will possess a frequency response like the ones in figure\,\ref{f:lumpedfresponse}, since the diaphragm acceleration cannot be expected to be constant throughout the frequency range.
Typically, a moving-coil dynamic transducer supplies the diaphragm motion, and the frequency response will then depend on the electrical and mechanical properties of the transducer and how they interact with the compression chamber.
The reason for providing the curves in figure\,\ref{f:lumpedfresponse} is to clarify the effects of the compression driver geometry; to obtain the actual frequency response of the system, a complete, coupled model is needed that also includes the transducer and a more realistic model for the waveguide load to replace assumption~\ref{ass}\,(v). 
\end{remark}

\section{Violations of assumptions~\ref{ass}}

In real compression drivers, assumptions~\ref{ass} will be violated to various degrees.
The sound pressure level in the chamber can be very high in the operational regimes common in professional audio, which is why the linearity assumption~(iv) can be questionable.
Although significant, this aspect is out of scope for the present contribution.
Another questionable postulate is that the diaphragm moves as a solid piston, assumption~(i).
In reality, structural modes form in the diaphragm at higher frequencies, which means that the effective moving mass will be frequency dependent, leading to an uneven frequency response.
This ``break-up'' mechanism may lower the moving mass at higher frequencies and therefore boost the high-frequency response, which, perhaps surprisingly, sometimes is regarded as beneficial. 
Also this (important) aspect will not be considered here.
However, violations of assumptions~(ii) and~(iii) will be addressed in detail.

\begin{figure}[b]\centering
\includegraphics{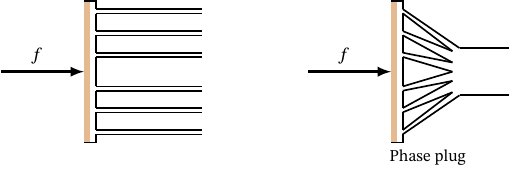}
\caption{Left: a conceptual cylindrical compression driver with annular openings, placed according to Smith's guidelines, to suppress radial modes in compression chamber.
Right: the function of the phase plug is to guide the sound from the slits into the throat of a horn.		
Cross-section view  through axis.}\label{FSmith}
\end{figure}

For real drivers, it does not hold that  the acoustic pressure $p$ is the same at each point within the chamber. 
In particular, the diameter of the diaphragm will typically be large enough so that acoustic modes form in the compression chamber.
In fact, the positioning of the outlet as in figure\,\ref{f:chamber} is particularly bad, mostly spoiling  the  effect of the compression chamber.
To address this issue, Smith~\cite{Sm53} developed an ingenious strategy that divides the outlet area $A_\text{wg}$ into $N$ circumferential slits.
If these are properly placed and sized, Smith showed that  the $N$ lowest radial modes in the chamber will not be excited.
An example with $N=3$ is shown to the left in figure\,\ref{FSmith}.
Such a configuration performs very closely to the ideal responses of figure\,\ref{f:lumpedfresponse} in the operational range below the first $N+1$ radial modes as long as the slits are perfectly absorbing.
However, in order to mount the driver to the mouth of a horn, which is typically cylindrical, the slits need somehow to be connected together using a device known as a \textit{phase plug}.
Figure\,\ref{FSmith} displays a simple design of a phase plug for $N=3$ Smith slits attached to the compression chamber.
Unfortunately, the phase plug introduces a nonideal load at the outlets from the compression chamber as well as  interdependencies between the loads to each slit, which means that the radial chamber modes will not be completely suppressed.
Moreover, the coupling between the slits means that internal resonances can also build up within the phase plug.
In a previous contribution~\cite{BeWaBe19}, we used gradient-based optimization to shape the phase-plug channels, which significantly improved the frequency response compared to a design as to the right in figure\,\ref{FSmith}, although we did not succeed in removing all resonances within the device. 

As we saw in \S\,\ref{Slumpedmodel}, high compression rates and narrow compression chambers are needed to obtain, for a given diaphragm movement, a high sound pressure and to avoid an excessive drop at high frequencies. 
The device will thus contain very narrow regions.
As a consequence, the viscothermal losses generated in the boundary layers cannot be ignored, violating the isentropic assumption.
This effect was noted already in the previous contribution~\cite{BeWaBe19}, where viscothermal losses were taken into account a posteriori by solving the linearized, compressible Navier--Stokes equations.
That is, the optimization was carried out using the isentropic Helmholtz equation for the acoustic pressure, but the  performance of the optimized device was assessed taking viscothermal losses into account.
Even though the final design outperformed the original design, this approach is conceptually unsatisfactory; it would be preferable to account for the losses 
also in the optimization step.

As discussed in the introduction, although the linearized, compressible Navier--Stokes equations accurately compute the lossy sound transmission, they are unpractical to use for shape optimization, due to extensive computational times and the need to manage the required extreme mesh refinements in the boundary-layer region.
Here we will instead rely on the boundary-condition model discussed in the introduction.
Berggren et al.~\cite{BeBeNo18} demonstrated the accuracy and extreme efficiency of this model on a compression driver with the same type of radial outlets as considered here.
The output sound pressure when employing the new model agreed closely with the one obtained from a more complete model, in which the linearized, compressible Navier--Stokes equations were used in the compression driver and phase plug.
However, the new model required about two orders of magnitude less memory and computational time compared to the more complete model.
Thus, the access to the boundary-layer model makes it now possible to apply numerical shape optimization to design the phase plug of a compression driver while taking viscothermal losses into account during the optimization.

\section{The compression driver with radial slits}

\begin{figure}[b]\centering
\begin{minipage}[t]{0.45\textwidth}\centering
\includegraphics[width=0.76\textwidth]{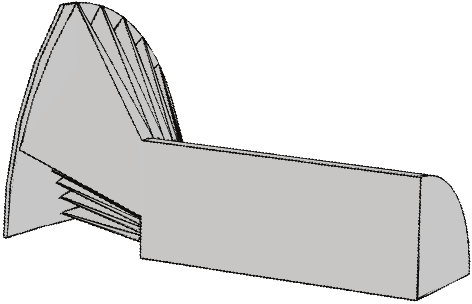}
\caption{Cutaway drawing of a compression driver with radial slits.}\label{Fradialslits}
\end{minipage}
\hfill
\begin{minipage}[t]{0.45\textwidth}
\includegraphics{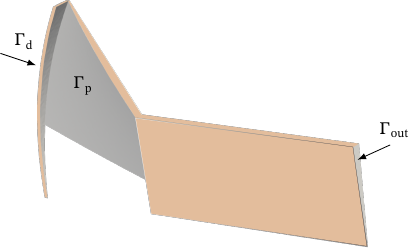}
\caption{The computational domain with the initial shape of the phase-plug wall $\Gamma_\text{p}$.
Symmetry conditions hold at $\Gamma_\text{sym}$, consisting of the full back side together with the surfaces indicated with tinted color in the figure.}\label{Finitdomain}
\end{minipage}
%
\end{figure}

An alternative to the Smith slits is to orient the outlets from the compression chamber in the  \textit{radial} direction, as in figure\,\ref{Fradialslits}.
Note that the wave propagation then is intrinsically three dimensional, whereas 2D axial symmetry modeling is possible  for  circumferential slits. 

The device we will consider contains $N_\text{s}$ radial slits, placed symmetrically in the circumferential direction, and the computational domain, depicted in figure\,\ref{Finitdomain}, contains one of these in half symmetry, with symmetry planes $\Gamma_\text{sym}$ placed in the radial direction.
The wall $\Gamma_\text{p}$ of the phase plug, shown in figure\,\ref{Finitdomain} in its initial shape, is the part of the device subject to shape optimization.
The phase plug connects the compression chamber, whose back side is the diaphragm $\Gamma_\text{d}$, to the waveguide.
On the spherically-shaped diaphragm $\Gamma_\text{d}$, we impose a uniform oscillatory movement  in the axial direction with angular frequency $\omega$.
This movement generates an acoustic pressure $P(\bs x, t) = \Real\e^{\i\omega t} p(\bs x)$ in the computational domain, and we require the complex pressure amplitude $p$ function to satisfy the following boundary-value problem for the Helmholtz equation,
\begin{subequations}\label{BVP}
\begin{alignat}{4}
-\Delta p - k^2 p &= 0		&&\text{in $\Omega$,}		\label{BulkHH}
\\
\i k p + \pder pn1 &= 0	&&\text{\,on $\Gamma_\text{out}$,}
\\
\pder pn1 &= 0									&&\text{on $\Gamma_\text{sym}$,}		\label{symbc}
\\
-\delta_V\frac{\i-1}{2} \DeltaT p + \delta_Tk^2\frac{(\i-1)(\gamma-1)}2 p + \pder{p}n1 &= 
\begin{cases}
\rho_0 a_\text{m}\,\bs n\cdot\bs e_\text{a}
\\
0
\end{cases}
&&\begin{aligned}
&\text{on $\Gamma_\text{d}$,}
\\
&\text{on $\Gamma_\text{w}\setminus\Gamma_\text{d}$,}
\end{aligned}\label{vtbc}
\\
\bs n_{\text T}\cdot\nablaT p&= 0		&&\text{on $\partial\Gamma_\text{w}$.}\label{BC2BC}
\end{alignat}
\end{subequations}
%
The boundary condition on the outlet $\Gamma_\text{out}$ of the waveguide is a simple radiation condition that perfectly absorbs planar modes.
In all numerical experiments, the frequency will be chosen so that all higher waveguide modes in the radial direction will be evanescent.
Following Oclee-Brown’s recommendations~\cite[Ch.\,7]{OB12}, we choose $N_s$ large enough to suppress circumferential modes in the compression chamber.
The symmetry boundary condition~\eqref{symbc} is imposed on the whole back side of the domain in figure\,\ref{Finitdomain} as well as on the tinted surfaces in the figure.
The purpose of the waveguide with the given boundary conditions is to impose a perfectly matched condition at the outlet of the phase plug; in reality, the outlet of the compression driver is typically directly mounted to the throat of a horn.

The viscothermal condition~\eqref{vtbc} is imposed on all remaining large solid wall surfaces $\Gamma_\text{w}$ of the domain: the diaphragm $\Gamma_\text{d}$ and the opposing surface inside the compression chamber, as well as the phase plug wall $\Gamma_\text{p}$.
Note, from figure\,\ref{Finitdomain}, that lossless homogeneous Neumann boundary conditions are imposed on the top part of the compression chamber as well as the top and bottom part of the phase plug, although viscothermal losses would be more appropriate.
This simplification streamlines the implementation and will only have a marginal acoustic effect.
In particular, the top and bottom parts of the phase plug are modeled without viscothermal losses in order to avoid cumbersome modifications to the shape gradient in the cases where the size of these parts changes.

The coefficients in condition~\eqref{vtbc} contain the viscous and thermal boundary-layer thicknesses
\begin{equation}\label{deltaVdeltaT}
\delta_V = \sqrt{\frac{2\nu}\omega}, \qquad
\delta_T = \sqrt{\frac{2\lambda}{\omega\rho_0 c_p}},
\end{equation}
where $\nu$ and $\lambda$ are the kinematic viscosity and thermal conductivity coefficient of air, respectively,  and $c_p$ the specific heat capacity of air at constant pressure.
This type of boundary condition, involving the \textit{tangential Laplacian operator} $\DeltaT$, can be viewed as a generalized impedance boundary condition, and is sometimes referred to as a \textit{Wentzell} boundary condition.
The tangential Laplacian operator in boundary condition~\eqref{vtbc} generates  diffusion in the tangential plane $\Gamma_{\text w}$ and can be defined as 
\begin{equation}
\DeltaT p = \nablaT\cdot\bigl(\nablaT p\bigr) 
\end{equation}
using the tangential gradient and tangential divergence operators
\begin{equation}\label{e:nablaT}
\nablaT p  = \nabla p - \bs n \pder pn1, \quad
\nablaT \cdot \bs q  = \nabla\cdot \bs q - \bs n \cdot \dfrac{\partial \bs q}{\partial n},
\end{equation}
defined for scalar and $\RR^3$-valued functions $p$ and $\bs q$, respectively.

Indeed, condition~\eqref{vtbc} constitutes in itself a diffusion--reaction problem on its surface, coupled to the bulk Helmholtz problem~\eqref{BulkHH} through the normal flux.
The presence of the tangential Laplacian requires a ``boundary condition to the boundary condition'', that is, a condition on the interface $\partial\Gamma_\text{w}$ towards neighboring boundary surfaces,  in this case $\Gamma_\text{sym}$.
We impose here condition~\eqref{BC2BC}, similarly as in the publication where the viscothermal model was devised~\cite{BeBeNo18}.
Vector $\bs n_\text{T}$ is the \textit{conormal} field at the boundary $\partial\Gamma_\text{w}$ of $\Gamma_\text{w}$.
(The conormal field is orthogonal to the normal field on $\Gamma_\text{w}$, outward-directed and in the tangent plane to $\Gamma_\text{w}$.)
The forcing given on $\Gamma_\text{d}$ in boundary condition~\eqref{vtbc}, in which $\bs e_\text{a}$ is a unit vector in the axial direction, stems from the diaphragm's oscillation with the acceleration amplitude $a_\text{m}$. 

The variational formulation of boundary-value problem~\eqref{BVP}, on which the finite-element discretization is based, is as follows.
\begin{equation}\label{varprob}
\begin{aligned}
&\text{Find $p\in W$ such that}
\\
&a(q, p) = l(q)\qquad\forall q\in W,
\end{aligned}
\end{equation} 
where
\begin{subequations}\label{adefldef}
\begin{align}
a(q, p) &= \int_\Omega\nabla q\cdot\nabla p - k^2\int_\Omega q p + \i k\int_{\Gamma_\text{out}} q p
+ \delta_Tk^2\frac{(\i-1)(\gamma-1)}2\int_{\Gamma_\text{w}} q p
+ \delta_V\frac{\i-1}2\int_{\Gamma_\text{w}}\nablaT q\cdot\nablaT p,
\label{adef}\\
l(q) &= \int_{\Gamma_\text{d}}\rho_0q a_\text{d} \,\bs n\cdot\bs e_\text{a},
\label{ldef}
\end{align}
\end{subequations}
where $W$, a strict subspace of $H^1(\Omega)$, is the closure of functions in $\mathscr C^1(\bar\Omega)$ with respect to the norm devised by
\begin{equation}\label{Wdef}
\lVert\ p \rVert^2_W = k_0^2\int_\Omega \lvert p\rvert^2 
+ \int_\Omega \lvert\nabla p\rvert^2 
+ \delta_T(\gamma-1)k_0^2\int_{\Gamma_\text{w}}\lvert p\rvert^2 
+ \delta_V\int_{\Gamma_\text{w}} \lvert\nablaT p\rvert^2.
\end{equation}

Note that condition~\eqref{BC2BC} is a ``natural condition'' for the tangential Laplacian, meaning that this interface will not enter into  the variational formulation, analogously to the vanishing of the symmetry condition~\eqref{symbc}.
Thus, this interface will be neither a source or sink of acoustic power.
The well-posedness proof in our previous article~\cite[Appendix~B]{BeBeNo18} also applies to problem~\eqref{varprob}.

\begin{remark}
To save space, we will in this article generally leave out the measure symbols in the integrals, like in expressions~\eqref{adefldef} and~\eqref{Wdef}.
The choice of integration measure will be clear from the domain of integration.
\end{remark}

\section{Optimization problem}

\begin{figure}[b]\centering
\includegraphics[width=0.3\textwidth]{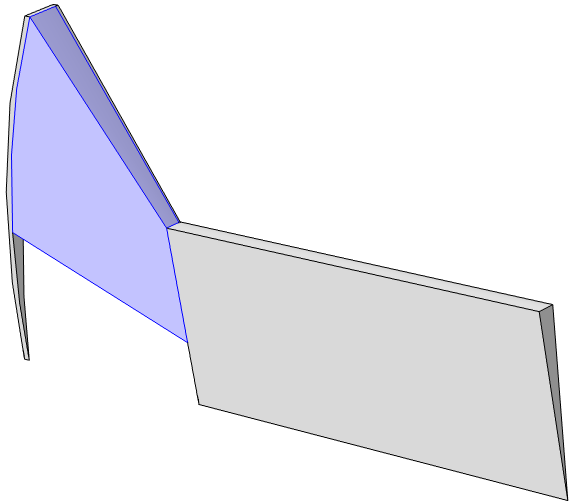}
\quad
\includegraphics[width=0.3\textwidth]{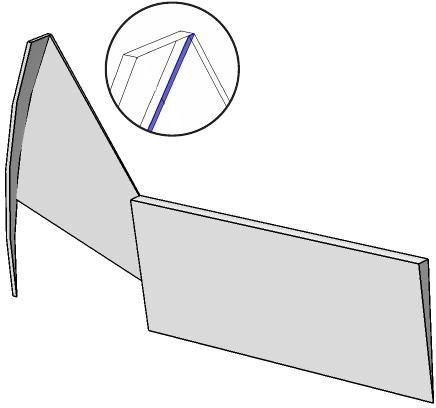}
\caption{The hold-all $\holdall$ is the union of the gray and tinted part  in the left picture.
The domain for the level set function $\phi$ is the tinted (blue) part $\Omega_\text{p, max}$.
The part of the computational domain that is fixed is the gray part, $\Omega_0 = \holdall\setminus\Omega_\text{p, max}$, shown to the right.
The tinted surface in the close-up indicates the exit from the compression chamber.
The computational domain is $\Omega=\Omega_0\cup\Omega_\text{p}$,  where $\Omega_\text{p}\subset\Omega_\text{p, max}$ such that  $\phi<0$.
}\label{FHoldall}
\end{figure}

The maximal domain $\holdall$ of all feasible computational domains,  the  \textit{hold all},  is depicted to the left in figure\,\ref{FHoldall}.
To describe the shape of the phase-plug boundary $\Gamma_\text{p}$, we  employ a level-set function  $\phi:\Omega_\text{p, max}\to\RR$, whose domain $\Omega_\text{p, max}\subset\holdall$ is the tinted part to the left in figure\,\ref{FHoldall}.
The complement $\Omega_0 = D\setminus\Omega_\text{p, max}$, illustrated to the right in figure\,\ref{FHoldall}, consists of the compression chamber,  the waveguide, and a thin region on the back side of the phase plug; $\Omega_0$ is fixed throughout the computations.

The phase-plug boundary and interior are defined through a level-set function on $\Omega_\text{p, max}$ as
\begin{equation}\label{GdOd}
\begin{aligned}
\Gamma_\text{p} &= \left\{\, \bs x\in\Omega_\text{p, max} \mid \phi(\bs x) = 0\,\right\},
\\
\Omega_\text{p} &= \left\{\, \bs x\in\Omega_\text{p, max} \mid \phi(\bs x) < 0\,\right\}.
\end{aligned}
\end{equation}
The level-set function will, in turn, implicitly be defined using a 3D generalization of the approach used in  previous studies~\cite{BeWaBe17,BeWaBe19}.
That is, $\phi$ is the solution to the boundary-value problem
\begin{subequations}\label{phiBVP}
\begin{alignat}{2}
-\Delta\phi &= \hat\phi			&&\qquad\text{in $\Omega_\text{p, max}$,} \label{phiBVPeq}
\\
\phi & = \phi_\text{D}					&&\qquad\text{on $\partial\Omega_\text{p, max}^\text{D}$,}		\label{PoissonDBC}
\\
\pder\phi n1 & = 0						&&\qquad\text{on $\partial\Omega_\text{p, max}\setminus\partial\Omega_\text{p, max}^\text{D}$.}		\label{PoissonNBC}
\end{alignat}
\end{subequations}
On the part of the boundary to $\Omega_\text{p, max}$ facing the compression chamber, we assign the Dirichlet boundary condition~\eqref{PoissonDBC}; the values of $\phi_\text{D}$ are given by the design illustrated in figure\,\ref{Finitdomain}, which also serves as the initial design for the optimization.
In this way, $\partial\Omega_\text{p, max}^\text{D}$ will always  be closed to the compression chamber, and the outlet from the compression chamber, shown as a tinted surface in the close up of the right picture in figure\,\ref{FHoldall}, will be a part of $\Omega_0$ and thus fixed throughout the computations.
We made this choice in order to fix a priori the compression ratio to be determined by other considerations not considered here.
(For instance, too high compression ratios and narrow chambers may introduce excessive sound pressures and high distortion.)
The homogeneous Neumann condition~\eqref{PoissonNBC} on the rest of the boundary yields freedom for the phase-plug boundary and interior to intersect the boundary of $\Omega_\text{p, max}$, a property that indeed will be observed in the optimized devices.
The level-set function will  depend on the right-hand-side function $\hat\phi$, whose nodal values in $\Omega_\text{p, max}$, after  finite-element discretization,  will be the actual design variables updated by the optimization algorithm.

The quantity of interest for the optimization is the mean sound pressure at the exit of the waveguide, that is,
\begin{equation}\label{e:qoi}
p_\text{out,k} = \frac1{\lvert\Gamma_\text{out}\rvert}\int\limits_{\Gamma_\text{out}} p,
\end{equation}
where $p$ is the solution to variational problem~\eqref{varprob} for wave number $k$, and $\lvert\Gamma_\text{out}\rvert$ the area of $\Gamma_\text{out}$.
Due to the absorbing boundary condition on $\Gamma_\text{out}$ and the evanescence of higher waveguide modes, the transmitted acoustic power is proportional to $\lvert p_{\text{out}, k}\rvert^2$.
For a set $\mathcal K$ of wave numbers, we will compare the results of optimization using the two objective functions
\begin{equation}\label{Jpow}
J_\text{pow}(\hat\phi) = \frac12\sum_{k\in\mathcal K}\frac1{\lvert p_{\text{out},k} \rvert^2}
\end{equation}
and 
\begin{equation}\label{Jtrack}
J_\text{track}(\hat\phi) = \frac12\sum_{k\in\mathcal  K} \bigl\lvert p_{\text{out},k}  - p^\text{ideal}_{\text{out}, k}\bigr\rvert^2,
\end{equation}
where 
\begin{equation}\label{poutidealdef}
p^\text{ideal}_{\text{out}, k} = p^\text{ideal}_k \e^{-\i kL} ,
\end{equation}
where $p^\text{ideal}_k$ is the ideal  pressure response as obtained from equation~\eqref{e:lumpedformula} and $e^{-\i kL}$ the phase factor associated with a planar wave propagating the distance $L$ between the diaphragm and the waveguide outlet $\Gamma_\text{out}$. 

Thus, by minimizing $J_\text{pow}$, we are maximizing the harmonic mean of the output power at the considered wave numbers, whereas by minimizing $J_\text{track}$, we are aiming to match the frequency response to the ideal one.
If needed, a Tikhonov regularization term may be added,
\begin{equation}\label{Jepsilon}
J_\epsilon(\hat\phi) = \epsilon J_\text{T}(\hat\phi) + J_\text{obj}(\hat\phi), 
\end{equation}
where $J_\text{obj}$ is either of the two objectives~\eqref{Jpow} or \eqref{Jtrack}, and 
\begin{equation}
J_\text{T}(\hat\phi) = \frac12 \int\limits_{\mathclap{\Omega_\text{d,max}}}(\hat\phi - \hat\phi_0)^2,
\end{equation}
in which $\hat\phi_0$ is the right-hand side function in equation~\eqref{phiBVPeq} associated with the starting guess of the optimization.
Such a regularization will put a bound on how much the level-set function can change from its initial shape.
In a previous study~\cite{BeWaBe17}, we found that this type of regularization was effective in limiting the curvature of the design boundary.
As we will see, in practice, it turns out that the Tikhonov regularization will not be necessary to employ in this case.

\section{Discretization}\label{s:discrete}

The cut finite element method (CutFEM)~\cite{BuClHaLaMa15}  is a discretization technique that is
similar to the XFEM approach~\cite{FriesBelytschko2000}.
These techniques are particularly suitable for problems, like here, in which the shape of the domain is unknown in the problem to be solved~\cite{DuysinxMiegroetJacobsEtAl2006,SharmaVillanuevaMaute2017,BurmanElfversonHansboEtAl2018,BretinChapelatOuttierEtAl2022}. 
We employ  a 3D version of the  approach  used in recent 2D acoustic shape optimization studies~\cite{BeWaBe17,BeWaBe19}.

The  hold all  $\holdall$ (figure\,\ref{FHoldall}) is triangulated into a mesh $\mathscr T_h$ of unstructured tetrahedra with maximum diameter $h>0$.
We require the mesh surfaces to conform with the interface between the level-set domain $\Omega_\text{p, max}$ and the rest of $\holdall$.
Let $V_h^n$ be the space of continuous functions that are polynomials of maximal degree $n$ on each tetrahedron $K\in\mathscr T_h$.
From $V_h^1$, we form the discrete space used to numerically solve the level-set boundary-value problem~\eqref{phiBVP}.
Let $U_h= V_h^1\big|_{\Omega_\text{p,max}}^{}$, and let $U_{h, 0}$ be the subspace of  functions  in $U_h$ vanishing on $\partial\Omega_\text{p,max}^\text{D}$.
The discrete level set function $\phi_h$ will then be the solution to 
\begin{equation}\label{phihprob}
\begin{aligned}
\text{$\phi_h\in U_h$ such that $\phi_h = \phi_\text{D}$ on $\partial\Omega_\text{p,max}^\text{D}$, and}
\\
\int\limits_{\mathclap{\Omega_\text{p,max}}}\nabla z_h\cdot\nabla\phi_h = \int\limits_{\mathclap{\Omega_\text{p,max}}} z_h\hat\phi_h
\qquad\forall z_h\in U_{h, 0},
\end{aligned}
\end{equation}
in which $\hat\phi_h\in\mathscr U_h$, where $\mathscr U_h$ is the space of design variables in the discrete case, which we choose to be equal to the space of test functions, $\mathscr U_h = U_{h, 0}$.

Then,  analogously as in definition~\eqref{GdOd}, the discrete level-set function $\phi_h$ defines the discrete phase-plug boundary and interior as
\begin{equation}\label{GdOdh}
\begin{aligned}
\Gamma_{\mathrm p, h} &= \left\{\, \bs x\in\Omega_\text{p, max} \mid \phi_h(\bs x) = 0\,\right\},
\\
\Omega_{\mathrm p, h} &= \left\{\, \bs x\in\Omega_\text{p, max} \mid \phi_h(\bs x) < 0\,\right\}.
\end{aligned}
\end{equation}
Note that $\Gamma_{\mathrm p, h}$ in general will cut through the inside of the elements of the triangulation, and  that $\Omega_{\mathrm p,h}$ thus will not conform to the mesh inside $\Omega_\text{p, max}$.
However, since $\phi_h$ is piecewise linear, the surface $\Gamma_{\mathrm p, h}$ will  be piecewise planar.
This property greatly simplifies both the assembly process for the acoustic finite-element matrices and the derivation of the shape calculus formulas for optimization.

\begin{figure}
\centering
\includegraphics[width=0.3\textwidth]{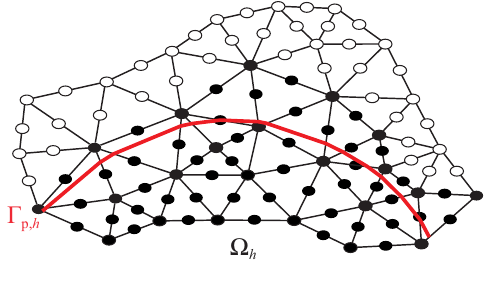}
\caption{The domain boundary $\Gamma_{\mathrm p, h}$ is defined by a vanishing  level set function $\phi_h$.
The computational domain $\Omega_h$ is here below $\Gamma_{\mathrm p, h}$. 
The nodes marked in black will affect the solution $p_h$ inside $\Omega_h$, whereas the white nodes will not. }\label{f:levelsetcut}
\end{figure}

The computational domain for the acoustic problem in the discrete case is $\Omega_h = \Omega_{\mathrm p, h}\cup\Omega_0$; $\Omega_{\mathrm p, h}$ will be a subset of the tinted volume ($\hat\Omega_\text{p,max}$) in figure\,\ref{FHoldall} and $\Omega_0$ is  the part shown to the right.
The finite-element approximation of variational problem~\eqref{varprob} will be carried out in the space
\begin{equation}\label{Whdef}
W_h = \left\{\, p_h \mid p_h = q_h|_{\Omega_h}^{} \text{for some $q_h\in V_h^2$}\,\right\}.
\end{equation}
Thus, we use piecewise-quadratic functions for the acoustic problem.
Note that the number of degrees of freedom for functions in $W_h$ will depend on the location of $\Gamma_{\mathrm p,h}$.
Also note that for elements that are intersected by $\Gamma_{\mathrm p,h}$, there will be basis function nodes for $W_h$ located outside of $\Omega_{\mathrm p,h}$.
This situation is illustrated in figure\,\ref{f:levelsetcut}, for simplicity of presentation in the analogous two-dimensional case.

The finite-element approximation of problem~\eqref{varprob} will be as follows.
\begin{equation}\label{acfem}
\begin{aligned}
&\text{Find $p_h\in W_h$ such that}
\\
&a(q_h, p_h) + \epsilon_\text{s} s_h(q_h, p_h) = l(q_h)\qquad\forall q_h\in W_h,
\end{aligned}
\end{equation} 
where $a(\cdot,\cdot)$ is defined as in expression~\eqref{adef} but with $\Omega = \Omega_h$. 
The bilinear term 
\begin{align}
\label{eq:ghost-penalty-def}
 s_h(p_h, q_h) = \sum_{S\in \mathscr{S}_h^g} \int_S 
 h^3 \llbracket \boldsymbol{n}  \cdot \nabla p_h \rrbracket
 \llbracket \boldsymbol{n}  \cdot \nabla q_h \rrbracket
\end{align}
is a so-called ghost penalty~\cite{Bu10,BuClHaLaMa15}, which is added to combat ill-conditioning of the system matrix in cases when only a small portion of an element intersected by $\Gamma_{\mathrm p,h}$ is located inside $\Omega_h$.
Here, $\epsilon_\text{s}> 0$ is a dimensionless stabilization parameter, $\mathscr{S}_h^g$ denotes the set of ghost penalty mesh faces, defined as the set of interior faces belonging to elements $K\in\mathscr T_h$ that are cut the embedded boundary $\Gamma_{p,h}$,
while $\llbracket \cdot \rrbracket$ denotes the jump of sufficiently smooth, element-wise defined functions across a face $S$.
More precisely, let $\bs f$ be a vector-valued function, satisfying $\bs f|_K\in C^0(\bar K)^3$ for each $K\in\mathscr T_h$.
Moreover, let $K_1, K_2\in\mathscr T_h$ be two elements sharing the mesh face $S$.
For any $\bs x\in S$ and for $i=1,2$, we define the limit functions
\begin{equation}
\bs f_i(\bs x) = \lim_{t\to0^+} \bs f(\bs x - t\bs n_i),
\end{equation}
where $\bs n_i$ is the outward-directed normal to $K_i$ on $S$.
Thus, $\bs f_1(\bs x)=\bs f_2(\bs x)$ only if $\bs f$ is continuous at $\bs x$.
Then the jump expression can be defined as
\begin{equation}\label{e:njump}
\llbracket \bs n\cdot\bs f(\bs x)\rrbracket = \bs n_1\cdot\bs f_1(\bs x) + \bs n_2\cdot\bs f_2(\bs x). 
\end{equation}
Since  $\bs n_1 = -\bs n_2$ on $S$, definition~\eqref{e:njump} indeed indicates a jump in the normal derivative of $\bs f$ at $S$.

The presence of ghost-penalty term $s_h$ does not affect the consistency of finite-element approximation~\eqref{acfem}, since it vanishes for sufficiently smooth $p$, where $p$ is the solution to variational problem~\eqref{varprob}.

\begin{remark}\label{Rextention}
Even though the number of active degrees of freedom in problem~\eqref{acfem} changes with the location of $\Gamma_{\mathrm p, h}$, it may be practical to keep the number of degrees fixed in the implementation.
A fixed number of degrees of freedom is easily accomplished by extending all functions in $W_h$ into $V_h$ by enlarging the system by an identity matrix and a zero right hand side for the degrees of freedom not affecting $p_h\large|_{\Omega_h}$; in figure\,\ref{f:levelsetcut} the solution would then vanish at the white nodes.
\end{remark}

Replacing $p_{\text{out}, k}$ in the definitions of objective functions~\eqref{Jpow}, \eqref{Jtrack}, and \eqref{Jepsilon} with corresponding numerical solution obtained from equation~\eqref{acfem}, we arrive at the following discrete optimization problem:
\begin{equation}
\begin{aligned}
&\text{Find $\hat\phi_h^*\in \mathscr U_h$ such that}
\\
&J_\epsilon(\hat\phi_h^*) \leq J_\epsilon (\hat\phi_h) \qquad\forall\hat\phi_h\in \mathscr U_h.
\end{aligned}
\end{equation}

\section{Shape calculus}\label{s:shapecalculus}

An evaluation of the function $\hat\phi_h\mapsto J_\text{obj}$ involves a composition of several operations,
\begin{equation}\label{compose}
\hat\phi_h \mapsto \phi_h \mapsto \Omega_h \mapsto p_h \mapsto p_{\text{out},k}\mapsto J_\text{obj}.
\end{equation}
That is, each function $\hat\phi_h$ gives rise to a level set function $\phi_h$ through the solution of problem~\eqref{phihprob}, which in turn through definition~\eqref{GdOdh} yields the computational domain $\Omega_h$ for problem~\eqref{acfem}.
Solutions to problem~\eqref{acfem} provides the acoustic pressures for each wavenumber, the quantities of interest $p_{\text{out}, k}$ in definition~\eqref{e:qoi},  and finally $J_\text{obj}$, either in terms of the power~\eqref{Jpow} or the tracking~\eqref{Jtrack} objective function.
To compute derivatives to be used by the optimization algorithm, we thus need to apply the chain rule to the composition~\eqref{compose}.

We have reported a shape calculus using the same kind of composition in earlier contributions~\cite{BeWaBe17, BeWaBe19}, except that the previous studies were done in 2D and without the boundary-layer model.
Apart from the mapping $\Omega_h\mapsto p_h$, previous calculations apply with minor adjustments.
However, the presence of the integrals over $\Gamma_{\mathrm w, h}$ arising from the boundary-layer model in bilinear form~\eqref{adef} significantly alters the shape calculus for this part of the mapping.
Tissot, Billard, and Garbard~\cite{TiBiGa20} carry out shape calculus for a conceptually similar case, that is, using a level-set geometry representation, cut elements, and the same viscothermal boundary-layer model.
However, their shape calculus relies on the classic formula~\eqref{e:dJ2formulaV} discussed in~appendix\,\ref{a:shapecalc}, which requires a smoothness that does not hold in the current case.
For instance, since the level-set function is continuous and piecewise linear, the design boundary will be continuous and piecewise planar but not $C^2$, and thus its curvature, which enters in formula~\eqref{e:dJ2formulaV}, can only be defined discretely as jumps of the normal field between element surfaces.
In addition, the smoothness assumption on the integrand, $f\in W^{1,1}(D)$, does not hold in our case either.

Instead, we will rely on a more general set of formulas discussed in~appendix\,\ref{a:shapecalc}.
Here we summarize the steps involved in computing the derivative of  the mapping  from the level-set function to the average acoustic pressure at the outlet,
\begin{equation}\label{e:j}
j(\phi_h) = p_{\text{out}, k}  = \frac1{\lvert\Gamma_\text{out}\rvert} \int_{\Gamma_\text{out}} p_h,
\end{equation}
and refer to~appendix\,\ref{a:shapecalc} for the derivation of the expressions.
Note that $p_{\text{out},k}$ is a complex number, and that the final objective function $J_\text{obj}$ (expression~\eqref{Jpow} or~\eqref{Jtrack}) is of least-squares type, involving a sum of squares over frequencies and over the real and imaginary parts.

The continuous, piecewise-linear level-set function $\phi_h$ can be expanded in a Lagrangian basis $\bigl\{\, w_i\,\bigr\}_{i=1}^M$.
Letting $w_l$, for $l\in\left\{\, 1, \ldots, M\,\right\}$,  be any of the basis functions of the expansion, we introduce the family of perturbations
\begin{equation}\label{e:phiht}
\phi_h^t = \phi_h + tw_l,
\end{equation}
parametrized by $t\geq0$. 
The perturbation generates a family of perturbed phase-plug wall boundaries and interiors
\begin{equation}\label{GhOht}
\begin{aligned}
\Gamma_{\mathrm p,h}^t &= \bigl\{\, \bs x\in\Omega_\text{p, max} \mid \phi_h^t(\bs x) = 0\,\bigr\},
\\
\Omega_{\mathrm p, h}^t &=  \bigl\{\, \bs x\in\Omega_\text{p, max} \mid \phi_h^t(\bs x) < 0\,\bigr\},
\end{aligned}
\end{equation} 
and the directional derivative of objective function~\eqref{e:j} is defined as
\begin{equation}\label{e:dj}
\textit dj(\phi_h; w_l) = \lim_{t\to0}\frac{j(\phi_h + tw_l) - j(\phi_h)}t.
\end{equation}

In most cases, we expect $\Gamma_{\text p, h}$ to cut through the interior of the mesh cells.
However, if $\phi_h$ happens to vanish along a mesh face, then this part of  $\Gamma_{\mathrm p, h}$ will coincide with that mesh face.
This is a rather special case, where, in fact, at most one-sided versions of directional derivative~\eqref{e:dj} can be expected to hold.
For simplicity, we exclude this case here as well as in the derivation in~appendix\,\ref{a:shapecalc}, an exclusion motivated by the fact that we did not experience any effects of this possible nonsmoothness in the numerical studies. 

Multiple times in the following, we will need to refer to the integrand for the terms in the variational form responsible for the modeling of viscothermal losses, that is, the terms in the integral over $\Gamma_\text{w}$ in definition~\eqref{adef}.
Thus, we will for convenience introduce the notation
\begin{equation}\label{e:psiqpndef}
\begin{aligned}
\psi(q, p; \bs n) &= \delta_Tk^2\frac{(\i-1)(\gamma-1)}2 q p + \delta_V\frac{\i-1}2\nablaT q\cdot\nablaT p
\\ &
= \delta_Tk^2\frac{(\i-1)(\gamma-1)}2  qp + \delta_V\frac{\i-1}2\bigl(\nabla q\cdot\nabla p - \bs n\cdot\nabla q\,\bs n\cdot\nabla p\bigr).
\end{aligned}
\end{equation}
Note that the dependency of $\psi$ on $\bs n$ is due to the definition~\eqref{e:nablaT} of the tangential gradient.

\begin{figure}\centering
\includegraphics{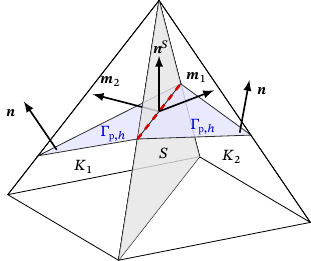}
\caption{Illustration of the geometrical objects involved in the evaluation  of the last term in the directional-derivative expression~\eqref{e:djcompute}.}\label{f:BdyIntersectS}
\end{figure}

As outlined in  \S\,\ref{s:discrete}, given a particular level-set function $\phi_h\in U_h$, the discrete acoustic pressure $p_h\in W_h$ is the solution of the finite-element problem
\begin{equation}
a(q_h, p_h) + \epsilon_h s_h(q_h, p_h) = l(q_h) 
\qquad\forall q_h\in W_h
\end{equation}
Letting $z_h\in W_h$ be the solution of the adjoint equation 
\begin{equation}
a(w_h, z_h) + \epsilon_h s_h(w_h, z_h) = \int_{\Gamma_\text{out}}\mkern-5mu w_h 
\qquad\forall w_h\in W_,
\end{equation}
the directional derivative~\eqref{e:dj} can be computed, with a notation explained below, from
\begin{equation}\label{e:djcompute}
\begin{aligned}
\dif j(\phi_h; w_l) &= \int\limits_{\Gamma_{\mathrm p, h}}\bigl(\nabla z_h\cdot\nabla p_h - k^2 z_h p_h\bigr)\frac{w_l}{\lvert\partial_n\phi_h\rvert}
\\
&\qquad+ \delta_V\frac{\i-1}2\int\limits_{\Gamma_{\mathrm p, h}}\Bigl(\frac{(P_T\nabla w_l)\cdot\nabla z_h}{\lvert\partial_n\phi_h\rvert}\pder{p_h}n1 + \pder{z_h}n1\frac{(P_T\nabla w_l)\cdot\nabla p_h}{\lvert\partial_n\phi_h\rvert}\Bigr)
\\
&\qquad+ \int\limits_{\Gamma_{\mathrm p, h}}\pder{}n1 \psi(z_h, p_h; \bs n)\, \frac{w_l}{\lvert\partial_n\phi_h\rvert}
+ \sum_{S\in\mathscr S_h} \int\limits_{\Gamma_{\mathrm p, h}\cap S} \bs n^S\cdot\bigl\llbracket \psi(z_h, p_h;\bs n) \bs m\bigr\rrbracket\frac{w_l}{\lvert\partial_{n^S}\phi_h\rvert},
\end{aligned}
\end{equation} 
where $P_T = I - \bs n\otimes\bs n$ is the projector on the tangent plane of $\Gamma_{\mathrm p, h}$, and where in the denominators, we use the shorthand notation $\partial_n$ for the normal derivative.
The last term in expression~\eqref{e:djcompute} is  of a sum over the set $\mathscr S_h$ of mesh faces, and the integral is evaluated over the  intersection between the boundary with the mesh surfaces.
Since we have assumed that no mesh surface coincides with a boundary portion, these intersections are all line segments.
Figure\,\ref{f:BdyIntersectS} illustrates a typical case, where it is to the left and right two mesh elements $K_1$, $K_2\in\mathscr T_h$ sharing the mesh face $S$, marked gray.
By assumption, the phase-plug boundary $\Gamma_{\mathrm p, h}$ passes through the interior of the mesh elements, and the piecewise-constant normal field $\bs n$ restricted on $\Gamma_{\mathrm p, h}\cap K_1$ and $\Gamma_{\mathrm p, h}\cap K_2$ is indicated in the figure.
The integration domain $\Gamma_{\mathrm p, h}\cap S$ is here the red dashed line segment.
The vector $\bs n^S$ is normal to $\Gamma_{\mathrm p, h}\cap S$  and in the plane $S$.
The conormal $\bs m_1$  is the vector in the plane $\Gamma_{\mathrm p, h}\cap K_1$ and outward-directed from $K_1$ and normal to $\Gamma_{\mathrm p, h}\cap S$; analogous conditions hold for $\bs m_2$. 
If $\Gamma_{\mathrm p, h}\cap K_1$ and $\Gamma_{\mathrm p, h}\cap K_2$  are coplanar,  then $\bs m_2 = -\bs m_1$.
The conormals can be used to define limits on $\Gamma_{\mathrm p, h}\cap S$ of functions on $\Gamma_{\mathrm p, h}$ with jump discontinuities over element boundaries.
That is, for $f:\Gamma_{\mathrm p, h}\to\RR$ such that, for $i=1, 2$,  $f|_{\Gamma_{\mathrm p, h}\cap K_i}\in C^0(\Gamma_{\mathrm p, h}\cap\bar K_i)$, we define, for each $\bs x\in\Gamma_{\mathrm p, h}\cap S$
\begin{equation}
f_i(\bs x) = \lim_{t\to 0^+} f(\bs x - t\bs m_i).
\end{equation}
Finally, the symbol $\llbracket \cdot\rrbracket$ in the last term of expression~\eqref{e:djcompute} denotes a jump of its argument over $\Gamma_{\mathrm p, h}\cap S$, defined by
\begin{equation}\label{e:mjump}
\bigl\llbracket \psi \bs m\bigr\rrbracket = \psi_1 \bs m_1 + \psi_2 \bs m_2.
\end{equation}
We note that the first term of $\psi$ in definition~\eqref{e:psiqpndef} will be continuous over $\Gamma_{\mathrm p, h}\cap S$ but not, in general, the second.

\begin{remark}
Note that we use, with a slight abuse of notation, the same symbol $\llbracket\cdot\rrbracket$ for different kinds of jumps over a face $S$ shared by two neighboring elements.
In definition~\eqref{e:njump}, the jump is with respect to limits associated with the normal to $S$, whereas in definition~\eqref{e:mjump}, the limits are associated with the conormals, at $S$, of the surface $\Gamma_{\text p, h}$ intersecting $S$.   
\end{remark}

\section{Test case specification and numerical results}

\begin{figure}
\hfill
\begin{minipage}{0.59\textwidth}
\includegraphics{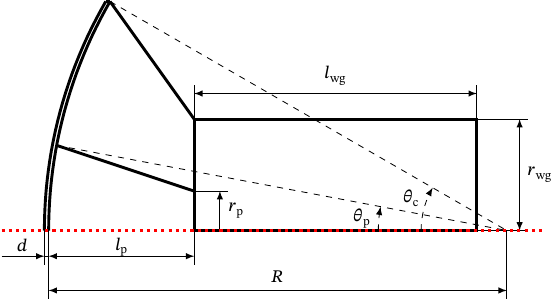}
\end{minipage}
\hfill
\begin{minipage}{0.3\textwidth}
\begin{tabular}{lllll}
\hline\hline
$\theta_\text{c}$	&$\pi/6$ & \quad\qquad&
$d$							& $0.5$ mm
\\
$r_\text{wg}$	& $13$ mm	&  &
$l_\text{wg}$	& $33$ mm 
\\
$l_\text{c}$					& $28$ mm & &
$l_\text{cp}$		&	$10$ mm	
\\
$r_\text{p}$		&	$4.6$ mm	& &$l_p$ & $17$ mm
\\
$R$ & $l_\text{c}/\theta_\text{c}$ & &
$\theta_\text{p}$ & $l_\text{cp}/R$
\\
\hline\hline
\end{tabular}
\end{minipage}
\hfill\null
\caption{Cross section of the maximal domain. }\label{Fcross-section}
\end{figure}

The compression driver we consider, illustrated in figure\,\ref{Fradialslits}, is equipped with $N_\text{s} = 32$ identical radial channels, and the computational domain is  illustrated in figure\,\ref{Finitdomain}.
The hold-all $D$, illustrated in figure\,\ref{FHoldall}, is constructed by rotating the planar domain in figure\,\ref{Fcross-section} the angle $\theta_{D} = \pi/32$ around the dotted axial  line.
The thin gray backside of the phase plug visible to the right in figure\,\ref{Fradialslits} is obtained by rotating corresponding planar part in figure\,\ref{Fcross-section} the angle $\theta_0=9.3773\cdot10^{-3}$ to fix the compression ratio to~$\kappa=12$.
The depth of the compression chamber is $d=0.5$\,mm.
The dimensions are chosen to be suitable as a driver for a horn operating in the upper midrange in the audio spectrum.
Table\,\ref{Tairprop} lists the air properties that are used in all computations to define the boundary-layer thicknesses~\eqref{deltaVdeltaT}.
The computations are carried out using 
the FEniCS computing platform~\cite{AlnaesBlechtaHakeEtAl2015} in combination with our in-house library extension \texttt{libCutFEM}, which implements the CutFEM-related algorithms described by Burman et al.~\cite{BurmanElfversonHansboEtAl2018}.

\begin{table}\centering
\caption{The air properties used in all computations.}\label{Tairprop}
\begin{tabular}{lll}
\hline\hline
Speed of sound & $c_0$	& $343.20$\,m s\textsuperscript{$-1$}	
\\
Air density & $\rho_0$	& $1.2044$\,kg m\textsuperscript{$-3$}
\\
Kinematic viscosity	& $\nu$ & $1.5061\times10^{-5}$\,m\textsuperscript2 s\textsuperscript{$-1$}	
\\
Prantdl number		& $N_\text{Pr}$ & $0.7078$
\\
Specific heat, constant pressure	& $c_p$ 	& $1.0049$\,kJ kg\textsuperscript{$-1$}K\textsuperscript{$-1$}
\\
Heat capacity ratio	& $\gamma$ & $1.4$
\\
Thermal conductivity & $\lambda$ & $c_p\mu/N_\text{Pr}$
\\
\hline\hline
\end{tabular}
\end{table}


\subsection{Baseline}

\begin{figure}\centering
\includegraphics{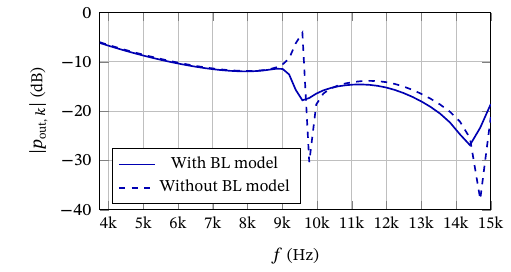}
%
\caption{Frequency response for the baseline design of figure\,\ref{Finitdomain}, with and without the boundary loss model.}\label{Fbaselinefres}
\end{figure}

For the baseline case, we fix a level set function $\phi_{h}$ that yields the phase plug design of figure\,\ref{Finitdomain}.
Here, the phase-plug wall $\Gamma_{\mathrm p, h}$ is a surface that attaches to the edge between the compression chamber and the thin phase plug back wall and, on the other side, to the front edge of the waveguide in figure\,\ref{FHoldall}. 
Using the air properties of table\,\ref{Tairprop} and holding the diaphragm oscillation amplitude $a_\textit{d}$ constant, we solve problem~\eqref{acfem} in the frequency range  $3.75$--$15$~kHz with as well as without the boundary loss model, the latter achieved simply by setting $\delta_T = \delta_V = 0$, in which case the model yields the standard hard-wall vanishing Neumann boundary condition, as is clear from expression~\eqref{vtbc}.
For a given, fixed value of the diaphragm acceleration $a_\text{d}$, figure\,\ref{Fbaselinefres} shows the frequency response in terms of the absolute value of the acoustic pressure at the outlet of the waveguide, $\lvert p_\text{out, $k$}\rvert$.
Two major resonances occur at $9$--$10$\,kHz and $14$--$15$\,kHz, respectively, both of which are significantly damped by viscothermal losses.
Thus, viscothermal effects cannot be neglected for frequencies higher than about 9\,kHz in this case.

\begin{remark}
Similarly as  already pointed out in remark~\ref{r:BiskopBrask1}, figure~\ref{Fbaselinefres} does not show a realistic frequency response of a full compression driver system, and neither will the response curves computed below for the optimized designs.
The calculation of realistic system frequency responses would require a coupled model including also the electromechanics of the transducer and a more realistic load at the waveguide. 
However, the purpose of the present study is to isolate the effects of the phase plug, not to simulate realistic frequency responses.
\end{remark}

\begin{figure}\centering
\includegraphics{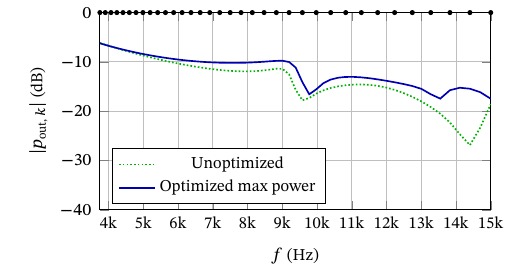}
\caption{Using objective function~\eqref{Jpow}: maximizing the harmonic mean of the output power at the 35~frequencies marked with black dots.}\label{Fmaxpow}
\end{figure}

\begin{figure}\centering
\includegraphics{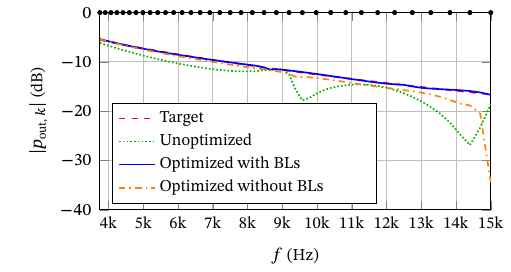}
	\caption{Using objective function~\eqref{Jtrack}: tracking the ideal frequency response at the 35~frequencies marked with black dots.}\label{Ftrack}
\end{figure}

\subsection{Optimization}

The optimization is carried out using the BFGS quasi-Newton algorithm as implemented in the SciPy Python library. 
The number of design variables, that is, the number of nodal values of the right-hand functions $\hat\phi_h$ in equation~\eqref{phihprob} is 4933.

Figure\,\ref{Fmaxpow} shows the frequency response of the optimized device when using  objective function~\eqref{Jpow}, that is, when maximizing the harmonic mean of the sound pressure $\lvert p_{\text{out}, k}\rvert$ at the waveguide outlet.
The mean is calculated over the 35 frequencies indicated with black dots in figure\,\ref{Fmaxpow} in the frequency range $[3.75, 15]$~kHz; each frequency is $2^{1/17}$ times the previous frequency.
The evaluation of the performance was done with the tighter frequency span of 69 frequencies, each pair spaced by the factor $2^{1/34}$.
The same frequency values for optimization and evaluation are used in all the subsequent cases below.
We note that the sound-pressure level is generally improved in the optimized device compared to the baseline.
However, the two resonances of the baseline are still clearly noticeable also in the optimized device.

This result indicates that a pure maximization of the sound pressure level may not be entirely satisfactory.
A better strategy could be to somehow give the optimization algorithm an incentive to reduce resonance and interference effects.
Therefore, it makes sense to consider the tracking approach, to minimize objective function~\eqref{Jtrack}.
Since the frequency response based on formula~\eqref{e:lumpedformula} is the best that can be achieved under the idealized condition of ignoring wave interference effects, we expect the optimization to counteract the interference effects in order to be able to track the idealized response. 
Figure\,\ref{Ftrack} shows the response of the baseline design in dotted green, the target frequency response $p_k^\text{ideal}$ from expression~\eqref{poutidealdef} in dashed red, and the response of the optimized device in solid blue using the tracking objective function~\eqref{Jtrack}.
The result is rather remarkable; the frequency response follows very closely the ideal one, without any noticeable resonance effects.

In the previous 2D study using radial phase plug slits~\cite{BeWaBe19}, the optimization was performed neglecting the viscothermal losses; these losses were accounted for only afterwards  when evaluating the performance of the optimized device.
In order to evaluate the necessity of the viscothermal model during optimization, we carried out the analogous strategy here.
We consider again the more successful tracking approach, optimizing using objective function~\eqref{Jtrack}.
The dash-dotted orange curve in figure\,\ref{Ftrack} shows the frequency response, evaluated using the viscothermal model but optimized without it. 
We note a significant deviation from the ideal response, particularly at the highest frequencies.

\begin{figure}
\begin{minipage}{0.49\textwidth}
\centering\includegraphics[width=0.4\textwidth]{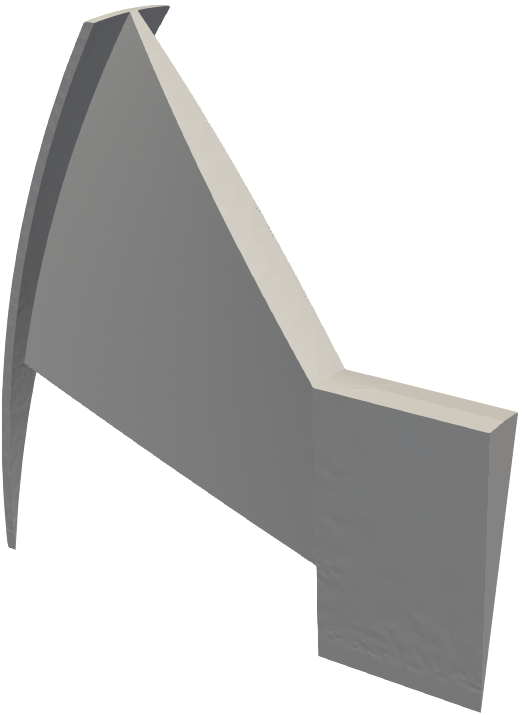}
\end{minipage}
\begin{minipage}{0.49\textwidth}
\centering\includegraphics[width=0.4\textwidth]{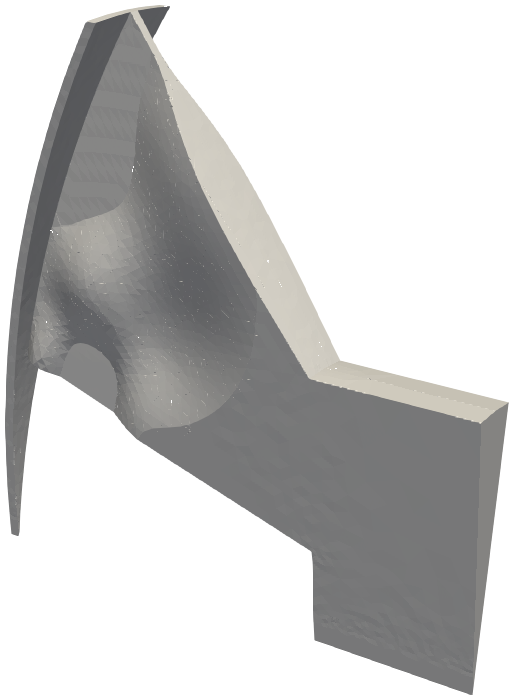}
\end{minipage}
\\[5mm]
\begin{minipage}{0.49\textwidth}
\centering\includegraphics[width=0.4\textwidth]{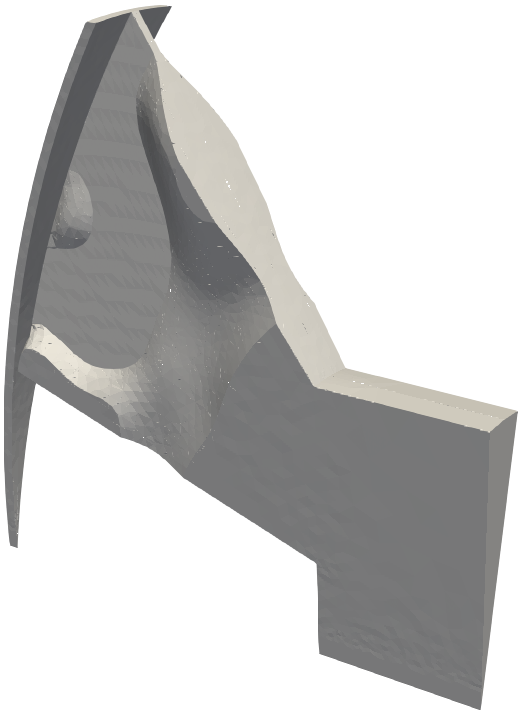}
\end{minipage}
\begin{minipage}{0.49\textwidth}
\centering\includegraphics[width=0.4\textwidth]{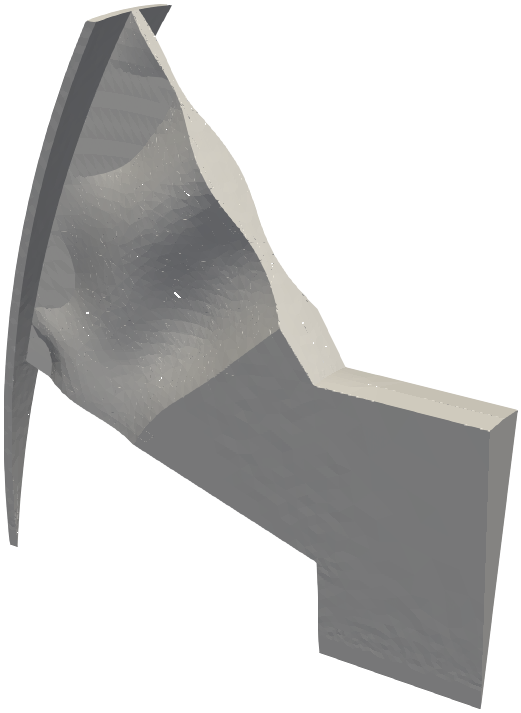}
\end{minipage}
\caption{Phase plug shapes.
Baseline (top left), optimized for maximum acoustic pressure at the waveguide outlet (top right), optimized by tracking the ideal response (bottom left), optimized by tracking the ideal response while ignoring the viscothermal losses during the optimization phase. 
For clarity, the waveguide has been shortened in the pictures.}\label{f:shapes}
\end{figure}

The shapes of all the considered phase plugs are visualized in figure\,\ref{f:shapes}, where the design is mirrored along the symmetry plane.
The use of the homogeneous Neumann condition~\eqref{PoissonNBC} in the boundary-value problem for the level-set function implies that the level-set function may intersect the boundary of $\Omega_\text{p, max}$. 
From figure\,\ref{f:shapes}, we see that this freedom is indeed utilized in all of the optimized designs.
In all of the optimized designs, note that neighboring phase plugs will be acoustically connected at the planar sections of the phase plug walls shown in figure\,\ref{f:shapes}.

\section{Discussion and conclusions}

The radial phase plug concept is  rarely used in today’s commercially available compression drivers, likely due to a lack of design guidelines  as straightforward as for annular channels.
Therefore, the use of numerical optimization techniques for shaping phase plugs is particularly interesting for this  case.
Perhaps surprisingly, when comparing with the results of shape optimization of annular channels~\cite{BeWaBe19}, it seems like it may be easier for the optimization to find an essentially resonance free design for radial channels, which may serve as an incentive for the industry to take a new interest in this arrangement.

Crucially important in obtaining the current results are the rather recent advances in basic modeling and optimization techniques.
For acoustic wave propagation in the narrow geometries inside a compression driver, it is necessary to account for viscothermal losses, also during the optimization, as shown by our reported results.
The recent popularization of the use of the Wentzell, or generalized impedance, condition to model viscothermal losses\,\cite{BeBeNo18}, is therefore a crucial component to be able to carry out the current study with a reasonable computational effort. 
Another advance utilized in this study is the combination of CutFEM methods with level-set geometry descriptions. 
In an earlier publication addressing the design of acoustic horns, we have shown that this approach gives exceptional design flexibility and may therefore find better designs than competing approaches~\cite{BeWaBe17}.
To obtain accurate gradient directions for rapid convergence of the optimization algorithm, we aim for an exact shape calculus for the discrete objective function, taking into account the lack of full regularity in the phase-plug wall geometry as well as the computed acoustic pressures. 
However, the fact that the viscothermal losses need to be modeled on the surfaces subject to design leads to complications in the shape calculus.
As far as we are aware, a fully discrete shape calculus involving the viscothermal boundary integral terms has not previously been carried out.
The analysis reported in~appendix\,\ref{a:shapecalc} uses a recent publication proving the relevant basic differentiation formulas~\cite{Be23}.



%
%
%

\section*{Acknowledgments}

Partial support for this work was provided by the Swedish Research Council, grant 2018-03546, and by eSSENCE, a strategic collaborative eScience program funded by the Swedish Research Council. 
The computations were performed on resources provided by the Swedish National Infrastructure for Computing (SNIC) at the High Performance Computing Center North (HPC2N).

\clearpage
\begin{appendices}
\section{Compression driver lumped-parameter model}\label{a:lumpedmodel}

We here provide a derivation of formula~\eqref{e:abslumpedformula} under assumptions~\ref{ass}.
The geometry of the compression driver is characterized by its rest volume $V_0$ and the areas $S_\text{d}$, $S_\text{wg}$ of its diaphragm and outlet. 
The easiest realization is the one of figure\,\ref{f:chamber}, although much more general geometries are admissible due to the lumped-parameter nature of the model.

Imagine a stiff massless membrane of area $S_\text{wg}$ positioned at the outlet of the chamber to the waveguide (the dashed tinted line in the cross-section view to the left in figure\,\ref{f:chamber}).
The velocity $U_\text{wg}$ of this membrane will then be the same as the air velocity at the outlet.
Under a movement of the diaphragm, the volume of the enclosed air between the piston and outlet membrane can be written as the sum of the rest volume $V_0$ plus an unsteady part $V'$,
\begin{equation}\label{Vexp}
V(t) = V_0 + V'(t),
\end{equation}
where the time derivative of the unsteady part satisfies
\begin{equation}\label{dotV'}
\dot V' = U_\text{wg}S_\text{wg} - U_\text{d} S_\text{d},
\end{equation}
in which $U_\text{d}$ and $U_\text{wg}$ are the velocities of the diaphragm and the massless membrane, respectively.
The values $U_\text{wg}(t)$ and $U_\text{d}(t)$ are  positive and negative, respectively, for velocities increasing the volume of the compression chamber.

By assumption~\ref{ass}(ii), the chamber will hold a uniform (but unsteady) pressure $P$, which satisfies, according to assumption~\ref{ass}(iii), the isentropy condition
\begin{equation}\label{isentropic}
PV^\gamma = \text{Const},
\end{equation}
where $\gamma$ is the heat capacity ratio.
Partitioning the pressure in the chamber into a steady and an unsteady part,
\begin{equation}\label{Pexp}
P(t) = p_0 + P'(t),
\end{equation}
substituting expansions~\eqref{Vexp}  and~\eqref{Pexp} into formula~\eqref{isentropic}, differentiating with respect to time, ignoring quadratic unsteady terms due to assumption~\ref{ass}(iv), and using expression~\eqref{dotV'}, we find that
\begin{equation}\label{linearized1}
\dot P' + \gamma\frac{p_0}{V_0}(U_\text{wg}S_\text{wg} - U_\text{d} S_\text{d}) =  0.
\end{equation}
Using that the speed of sound in an ideal gas (assumption~\ref{ass}(iii)) satisfies
\begin{equation}
c_0 ^2  = \gamma \frac{p_0}{\rho_0},
\end{equation}
and definitions
\begin{equation}
d = \frac{V_0}{S_\text d}, \qquad \kappa = \frac{S_\text d}{S_\text{wg}},
\end{equation}
where $d$ is the depth of the chamber in the simple geometry of figure\,\ref{f:chamber} and $\kappa$ the compression ratio, expression~\eqref{linearized1} can be written

\begin{equation}\label{linearized2}
\dot P' + \frac{\rho_0 c_0^2}d \left(\frac{U_\text{wg}}\kappa - U_\text d\right) = 0.
\end{equation}
Under time-harmonic excitation with frequency $\omega$ (assumption~\ref{ass}(i)),
\begin{equation}
P'(t) = \Real\bigl(p \e^{\i\omega t}\bigr), \qquad
U_\text{wg}(t) = \Real\bigl(u_\text{wg} \e^{\i\omega t}\bigr), \qquad
U_\text{d}(t) = \Real\bigl(u_\text{d}\e^{\i\omega t}\bigr), 
\end{equation}
expression~\eqref{linearized2} becomes
\begin{equation}\label{linearized3}
\i\omega p + \frac{\rho_0 c_0^2}{d\kappa} u_\text{wg} = \frac{\rho_0 c_0^2}d u_\text{d}
= \frac{\rho_0c_0^2}{d\,\i\omega}a_d,
\end{equation}
replacing velocity with acceleration amplitude in the last equality.
Due to the perfect termination of the waveguide (assumption~\ref{ass}(v)),  the acoustic pressure and velocity in the waveguide are proportional with proportionality constant the specific characteristic impedance $\rho_0 c_0$.
Substituting $p = \rho_0 c_0 u_\text{wg}$ into equation~\eqref{linearized3} and multiplying with $\i\omega d/c_0^2$, we find the requested formula 
\begin{equation}
\left(-d k^2 + \i\frac k\kappa\right)p = \rho_0 a_d,
\end{equation}
where $k = \omega/c_0$ is the wave number.

\section{Shape calculus details}\label{a:shapecalc}

Here we detail the calculations that yields expression~\eqref{e:djcompute}, that is, the directional derivative of functional~\eqref{e:j} with respect to the perturbation~\eqref{e:phiht} of the level-set function.

\subsection{Basic formulas for volume and surface integrals}\label{s:formulas}

The basic toolbox of shape calculus contains expressions of directional derivatives of integrals over a domain $\Omega\subset\RR^d$ or its boundary $\partial\Omega$. 
The directional derivatives are defined with respect to \textit{domain paths} $t\to\Omega_t$, where $\Omega_t$ is a family of domains parametrized by a scalar $t$ in which $\Omega_0 =\Omega$.
In the literature, such as in the classic monographs by Delfour \& Zolésio~\cite{DeZo11} and Sokolowski \& Zolésio~\cite{SoZo92}, shape calculus is typically carried out assuming that the domain path is generated by a smooth homeomorphism $T_t:D\to\RR^d$, where $D\subset\RR^d$ is a \textit{hold-all} containing all admissible domains.
Thus, $\Omega_t = T_t(\Omega)$, and using a change of variables based on $T_t$, integrals over $\Omega_t$ and $\partial\Omega_t$ can be converted into integrals over $\Omega$ and $\partial\Omega$.
The transformation parameters then appears explicitly in the integrands, which therefore can be differentiated using ordinary calculus.
A particularly popular choice to generate the transformation is through a \textit{perturbation of identity}, so that for any $\bs x\in\Omega$, $T_t(\bs x) = \bs x + t\bs V(\bs x)$, where $\bs V: D\to\RR^d$ is a vector field.
The use of this type of transformation marries well with the need to modify a computational mesh to changes in the domain boundary.
The so-called velocity field $\bs V$ can then be used to displace the boundary mesh nodes as well as the internal nodes, the latter to preserve mesh quality.

The use of transformations to generate domain paths is less natural, however, in the current context. 
All admissible phase-plug shapes are confined within the hold-all $D$ in figure\,\ref{FHoldall}, but the computational mesh is fixed, the phase-plug wall is a boundary internal to $D$, and there are no transformations involved in the movement of the phase-plug wall. 
Nevertheless, from perturbation~\eqref{e:phiht}, a (nonunique) artificial velocity field $\bs V$ can be defined locally on the mesh elements that are cut by the zero set of the level-set function (cf.~figure\,\ref{f:levelsetcut}).
The shape calculus carried out in our previous work using a 2D CutFEM approach~\cite{BeWaBe19,BeWaBe17} relied on such an approach. 

However, it is not necessary to invoke transformations and artificial velocity fields in this case.
An alternative is to generalize a method introduced by Delfour~\cite{De18} and consider \textit{dilations}, generated by perturbation~\eqref{e:phiht}, of surface patches defined by the level-set function.
We will here rely on formulas for the shape derivative of domain and boundary integrals as derived, using such a method,  in a recent contribution by the first author~\cite{Be23}.
These formulas assume minimal smoothness and takes into account that the design boundary will contain edges and a discontinuous normal field when the design boundary intersects mesh surfaces, as illustrated in figure\,\ref{f:BdyIntersectS}.
The presence of such discontinuities did not affect the shape calculus in our previous work, due to the natural boundary condition assumed there.
Here, however, the viscothermal boundary terms in the bilinear form~\eqref{adef}  contains discontinuities that need to be handled.

\begin{figure}\centering
\begin{minipage}[t]{0.4\textwidth}\centering
\includegraphics{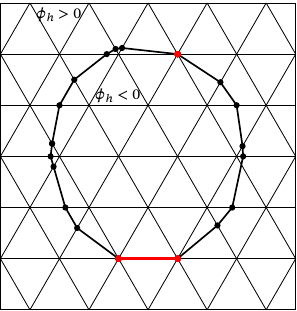}
\end{minipage}
\begin{minipage}[t]{0.4\textwidth}\centering
\includegraphics{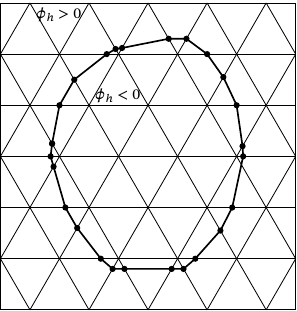}
\end{minipage}
\caption{Examples in 2D of a meshed rectangular hold-all $D$ and  domains $\Omega_h$ generated by  level-set functions $\phi_h$ defined on the mesh in the figure.
Both types of domains are admissible for theorem~\ref{t:dJ1wthm}, although only one-sided derivatives exists for domains like the one to the left, where a mesh face is contained in the boundary.
However, the proof~\cite[Theorem~6.11]{Be23} of theorem~\ref{t:dJ2wthm} requires domains like the one to the right, where the domain boundary does not intersect any mesh nodes.}\label{FDomain}
\end{figure}

We thus consider a continuous, piecewise-linear level-set function $\phi_h$, defined on a fixed triangulation $\mathscr T_h$ on $D$,  generating the computational domain $\Omega_h$, defined by the condition $\phi_h<0$. 
Consider then a perturbation
\begin{equation}\label{e:phit}
\phi_{h,t} = \phi_h + t w,
\end{equation}
where $w$ is a Lagrangian basis function and $t>0$.
This perturbation generates a perturbed domain $\Omega_{h,t}\subset\Omega_h$; the inclusion is due to that each basis function $w$ for continuous, piecewise-linear functions satisfies $w\geq0$.
The following formula for the shape derivative of volume integrals over $\Omega_h$ holds.
\renewcommand{\thesection}{\Alph{section}}
\begin{theorem}\label{t:dJ1wthm}
Under perturbation~\eqref{e:phit} and for $t\mapsto f(t)$ and $t\mapsto f'(t)$ continuous in some nonempty interval $[0, t_\text{max}]$ such that $f(t), f'(t) \in C^0(\bar{\mathscr T}_h)$ on $(0,  t_\text{max})$, the directional semiderivative of volume integral 
\begin{equation}
J_1(\phi_{h,t}) = \int_{\Omega_{\mathrlap{h,t}}} f_t\,\textit{dV}
\end{equation}
at $t=0$ satisfies
\begin{equation}\label{e:dJ1formula}
\textit dJ_1(\phi_h; w) = 
\lim_{t\to0^+} \frac1t\bigl( J_1(\phi_{h,t})  - J_1(\phi_h)\bigr) =
\int_{\Omega_{\mathrlap h}} f'\,\textit{dV} - \int_{\partial\Omega_h} f\frac{w}{\lvert\partial_n\phi_h\rvert}\,\textit{dS}.
\end{equation}
\end{theorem}
\begin{remark}
The notation $g\in C^k(\bar{\mathscr T_h})$ means that $g|_K\in C^k(\bar K)$ for each element $K\in\mathscr T_h$.
\end{remark}
Theorem~\ref{t:dJ1wthm} is an immediate consequence of the same theorem proven for a fixed $f\in C^0(\bar{\mathscr T}_h)$~\cite[Theorem~6.7]{Be23}.
Note that the theorem only refers to a one-sided derivative; the limits  $t\to0^+$ and $t\to0^-$ both exist but may be different.
The reason for this lack of strong differentiability is that $f$ as well as $\partial_n\phi_h$ are allowed to contain jump discontinuities across mesh faces.
These discontinuities become an issue for domains like the one to the left in figure\,\ref{FDomain},  where $\int_{S\cap\partial\Omega}\,\textit{dS}>0$ for some mesh face $S$.
For basis functions $w$ associated with the end points of such an $S$, formula~\eqref{e:dJ1formula} holds for the limit of $f$ and $\partial_n\phi_h$ on $S$ from the \textit{interior} of $\Omega_h$;
(Recall that $\Omega_{h,t}\subset\Omega_h$ for $t\geq0$.)
When the limit $t\to0^-$ is considered, it will be the limit of the values from the exterior that should be employed, since $\Omega_{h,t}\supset\Omega_h$ for $t\leq0$.
However, when $\partial\Omega_h$ only cuts through the interior of the mesh elements, as in the right picture of figure\,\ref{FDomain}, the limits as $t\to0^+$ and $t\to0^-$ will agree, since then $f$ and $\partial_n\phi_h$ will vary continuously under small enough perturbations of the boundary.

To carry out the shape calculus, we will need, in addition to formula~\eqref{e:dJ1formula}, also the corresponding formula for integrals over $\partial\Omega_h$.
Note that any nonempty boundary segment $\partial\Omega_h\cap K$, for some tetrahedral element $K\in\mathscr T_h$, will either be a triangle or a quadrilateral. 
At the edges of each such $\partial\Omega_h\cap K$, we will need to use the concept of conormals, introduced in §\,\ref{s:shapecalculus}.
These will be used to define limits and jumps of functions across line segments $\partial\Omega_h\cap S$, where $S$ is a mesh face, as  illustrated in figure\,\ref{f:BdyIntersectS}.

The following theorem is an immediate consequence of the same theorem proven for a fixed $f\in C^1(\bar{\mathscr T_h})$~\cite[Theorem~6.11]{Be23}.
Similarly as for theorem~\ref{t:dJ1wthm}, nondifferentiability can be expected when  some mesh face $S$ aligns with $\partial\Omega_h$, like in the domain to the left in figure\,\ref{FDomain}.
However, the ambiguities due to the discontinuities in the integrand of the last term in the expression for $\textit dJ_2$ below are not as easily resolved as in theorem~\ref{t:dJ1wthm}.
 Thus,  the proof of theorem~\ref{t:dJ2wthm} (\cite[Theorem~6.11]{Be23}) is restricted to domains where the boundary does not intersect any mesh nodes, such as the one to the right in figure\,\ref{FDomain}.
\begin{theorem}\label{t:dJ2wthm}
Assume that the boundary $\partial\Omega_h$ of the domain $\Omega_h$ does not intersect with any mesh nodes. 
Then,  under perturbation~\eqref{e:phit} and for $t\mapsto f(t)$ and $t\mapsto f'(t)$ being continuous in some nonempty interval $[0, t_\text{max}]$ such that $f(t), f'(t) \in C^1(\bar{\mathscr T}_h)$ on $(0,  t_\text{max})$, the directional derivative of surface integral 
\begin{equation}
J_2(\phi_{h,t}) = \int_{\partial\Omega_{\mathrlap{h,t}}} f_t\,\textit{dS}
\end{equation}
at $t=0$ satisfies
\begin{equation}\label{e:dJ2formula}
\textit dJ_2 (\phi_h, w) =  \int_{\partial\Omega_h}\Bigl( f'-\pder fn1 \frac w{\lvert\partial_n\phi_h\rvert}\Bigr)\,\textit{dS}
- \sum_{S\in\mathscr S_h}\int_{\partial\Omega_h\cap S} \bs n^S\cdot \llbracket f\bs m
\rrbracket  \frac w{\lvert\partial_{n^{S}}\,\phi_h\rvert}\,\textit{d}\gamma,
\end{equation}
where $\bs n^S$ is the normal vector to $\partial\Omega_h\cap S$, located in $S$ and outward-directed from $\Omega_h$,  
and where $\llbracket f\bs m\rrbracket = f_1\bs m_1 + f_2\bs m_2$.
Here, for $k= 1, 2$, $\bs m_k$ are the conormals to $\partial\Omega\cap K_k$ at $\partial\Omega\cap S$, where $K_1$, $K_2\in\mathscr T_h$ such that $\bar S = \bar K_1\cap\bar K_2$.
Each $\bs m_k$ lies in the plane of and is directed outward from $\partial\Omega_h\cap K_k$, and  $f_k$ is the limit $f$ on $S$ defined by $f_k(\bs x) = \lim_{\epsilon\to0^+} f(\bs x - \epsilon\bs m_k)$ , for any $\bs x\in S$.
\end{theorem}

It is  illuminating to compare directional derivatives~\eqref{e:dJ1formula} and~\eqref{e:dJ2formula} with the traditional expressions, that is, the formulas obtained when utilizing domain transformations through a perturbation of identity. 
For domain integrals, the classical formula is
\begin{equation}\label{e:dJ1formulaV}
\textit J_1(\Omega; \bs V) = \int_\Omega f' \,\textit dV +\int_{\partial\Omega} f \bs V\cdot\bs n\, \textit dS
\end{equation}
The proof by Delfour \& Zolésio~\cite[Ch.\,9, Theorem~4.2]{DeZo11} holds for Lipschitz domains $\Omega$, velocity fields satisfying $\bs V\in C^1(D)^d$ and for functions $t\mapsto f_t$ continuous with values in $W^{1,1}(D)$ and differentiable with values in $L^1(D)$.
For boundary integral, the classical formula is
\begin{equation}\label{e:dJ2formulaV}
\textit J_2(\Omega; \bs V) = \int_{\partial\Omega} f'\,\textit dS + \int_{\partial\Omega}\biggl(\pder fn1 + \kappa f\biggr)\bs V\cdot\bs n\,\textit dS,
\end{equation}
where $\kappa$ is the local summed curvature (twice the mean curvature) of the boundary.
The conditions for the boundary integral formula is quite restrictive.
The proof by Delfour \& Zolésio~\cite[Ch.\,9, Theorem~4.2]{DeZo11}  holds for $\Omega$ of class $C^2$ and for continuous  functions $t\mapsto f_t$ with values in $H^2(D)$. 

Comparing our formulas~\eqref{e:dJ1formula} and~\eqref{e:dJ2formula} with the classical~\eqref{e:dJ1formulaV} and~\eqref{e:dJ2formulaV}, we first note that the term~$-w/\lvert\partial_n\phi_h\rvert$ corresponds to the normal velocity $\bs V\cdot\bs n$ in the classical formulas.
Secondly, we note that the last integral in expression~\eqref{e:dJ2formula} is a discrete counterpart to the curvature term in expression~\eqref{e:dJ2formulaV}.
Note that the requirement for formula~\eqref{e:dJ2formulaV} of a $C^2$-regular domain and and integrand $f\in H^2(D)$ is violated in our case.

\subsection{Calculation of directional derivative~\eqref{e:djcompute}}

Based on basic formulas of theorems~\eqref{t:dJ1wthm} and~\eqref{t:dJ2wthm}, we here provide the details of the calculation of directional derivative~\eqref{e:djcompute}.
Due to the differentiability subtleties associated with theorem~\ref{t:dJ2wthm}, we will for the calculations below assume that the boundary $\partial\Omega_h$ does not intersect any mesh nodes, as in the domain to the right in figure\,\ref{FDomain}.

Recall from \S\,\ref{s:discrete} that the solution space $W_h$ to problem~\eqref{acfem} will not be fixed throughout the optimization, since it is spanned by the basis functions whose support intersect  domain $\Omega_h$ with positive measure.
As illustrated in  figure\,\ref{f:levelsetcut}, the nodal degrees of freedom for $W_h$, marked black in the figure, are located not just inside $\Omega_h$ but also outside the domain for the elements whose interior intersects the boundary.
Considering now perturbation~\eqref{e:phiht} for $t\geq0$, we  note the nodal degrees of freedom for $W_h$ only changes  when the domain boundary moves away from an intersection with a mesh node. 
Thus, under the assumption made above — that the boundary $\partial\Omega_h$ does not intersect any mesh node — we conclude that when starting with $t$ sufficiently small,  $W_h$ will not change during the limit process $t\to0$.

Thus, given the solution space $W_h$, there is a $t_\text{max}>0$  such that for each $t\in[0, t_\text{max}]$, the discrete acoustic pressure on the perturbed domain will satisfy the following problem.
\begin{equation}\label{e:acfemperturb}
\begin{aligned}
&\text{Find $p_h^t\in W_h$ such that $\forall q_h\in W_h$,}
\\
&\int_{\Omega_h^t}\!\!\nabla q_h\cdot\nabla p_h^t - k^2\int_{\Omega_h^t}\!\! q_h p_h^t + \i k\int_{\Gamma_{\mathrlap{\text{out}}}} q_h p_h^t
+ \int_{\Gamma_{\text w, h}^t}\!\!\!\psi(q_h,p_h^t; \bs n^t) + \epsilon_h s_h(q_h,p_h^t) = \int_{\Gamma_{\mathrlap{\text{d}}}}\rho_0q_h a_\text{d} \,\bs n\cdot\bs e_\text{a}.
\end{aligned}
\end{equation} 
where $\psi$ is defined in expression~\eqref{e:psiqpndef}.
Consequently, the average acoustic pressure at the outlet becomes
\begin{equation}\label{e:jphit}
j(\phi_h^t)  = \frac1{\lvert\Gamma_\text{out}\rvert} \int_{\Gamma_{\mathrlap{\text{out}}}}p_h^t.
\end{equation}
Note that outlet boundary $\Gamma_\text{out}$ is unchanged by the perturbation.
The \textit{shape derivative} of $p_h^t$ at $t=0$, defined by
\begin{equation}\label{shapeder}
p_h' = \lim_{t\to0} \frac1t\bigl(p_h^t - p_h^0\bigr),
\end{equation}
is also an element in $W_h$, since $p_h^t\in W_h$ for each $t\in[0, t_\text{max}]$, 
Differentiating objective function~\eqref{e:jphit} at $t=0$ yields
\begin{equation}\label{e:djp'}
\textit dj(\phi_h; w_l) = \int_{\Gamma_{\mathrlap{\text{out}}}}p_h'.
\end{equation}

To determine an explicit expressions for directional derivative $\textit dj$ in terms of $w_l$, we start by differentiating equation~\eqref{e:acfemperturb} at $t=0$, utilizing formulas~\eqref{e:dJ1formula} and~\eqref{e:dJ2formula},  to find that, for each $q_h\in W_h$,
\begin{equation}\label{e:acfemdiff}
\begin{aligned}
&\int_{\Omega_{\mathrlap h}} \nabla q_h\cdot\nabla p_h' - k^2\int_{\Omega_{\mathrlap h}} q_h p_h' 
- \int_{\partial\Omega_{\mathrlap{h}}}\bigl(\nabla q_h\cdot\nabla p_h - k^2 q_h p_h\bigr)\frac{w_l}{\lvert\partial_{\!n}\,\phi_h\rvert}
+ \i k\int_{\Gamma_{\mathrlap{\text{out}}}} q_h p_h'
\\
&\qquad
+ \int_{\Gamma_{\mathrlap{\mathrm w, h}}} \psi'(q_h,p_h; \bs n)   - \int_{\Gamma_{\mathrm w, h}}\! \pder{}n1\psi(q_h,p_h;\bs n) \frac{w_l}{\lvert\partial_{\!n}\,\phi_h\rvert}
- \sum_{S\in\mathscr S_h} \int\limits_{\Gamma_{\mathrm w, h}\cap S} \bs n^S\cdot\bigl\llbracket \psi(q_h, p_h;\bs n) \bs m\bigr\rrbracket\frac{w_l}{\lvert\partial_{n^{\!S}}\phi_h\rvert} 
\\
&\qquad\qquad
+ \epsilon_h s_h(q_h, p_h')= 0.
\end{aligned}
\end{equation}

Here, similarly as previously observed in the analogous 2D case~\cite{BeWaBe19}, we have used that the mapping $t\mapsto s_h(q_h, p_h^t)$ is differentiable, which is a consequence of our assumption that the boundary $\partial\Omega_h^t$ does not intersect any mesh node. 
Differentiating $\psi$, using definition~\eqref{e:psiqpndef}, we find that
\begin{equation}\label{e:phi'}
\begin{aligned}
\psi'(q_h, p_h; \bs n)  &= \der{}t1\psi(q_h,p_h^t; \bs n^t)\big|_{t=0}
\\
&= \alpha_T q _hp'_h + \alpha_V\bigl(\nabla q_h\cdot\nabla p'_h - \bs n'\cdot\nabla q_h\pder {p_h}n1 - \pder {q_h}n1\bs n'\cdot\nabla p_h - \pder {q_h}n1\pder {p'_h}n1\bigr)
\\ &
= \alpha_T q_h p'_h + \alpha_V\nablaT q_h\cdot\nablaT p'_h 
- \alpha_V\bigl(\bs n'\cdot\nabla q_h\pder{p_h}n1 + \pder{q_h}n1\bs n'\cdot\nabla p_h\bigr).
\end{aligned}
\end{equation}
Notice that $\bs n'$ vanishes everywhere  except on the part of the phase-plug boundary affected by the perturbation, $\Gamma_{\mathrm p, h}$ , on which the perturbed normal field satisfies
\begin{equation}\label{e:nt}
\bs n^t = \frac{\nabla(\phi_h + t w_l)}{\bigl\lvert\nabla(\phi_h + t w_l)\bigr\rvert}\qquad \text{on $\Gamma_{\mathrm p, h}^t$}.
\end{equation}
We will differentiate the right side of expression~\eqref{e:nt}, starting with the denominator
\begin{equation}
\begin{aligned}
\der{}t1 \bigl\lvert\nabla(\phi_h + t w_l)\bigr\rvert\Big|_{t=0}
&= \der{}t1 \bigl(\nabla(\phi_h + t w_l)\cdot\nabla(\phi_h + t w_l)\bigr)^{1/2}\Big|_{t=0}
\\ &
= \bigl(\nabla(\phi_h + t w_l)\cdot\nabla(\phi_h + t w_l)\bigr)^{-1/2}\nabla(\phi_h+tw_l)\cdot\nabla w_l\Big|_{t=0}
\\ &
=\frac{\nabla\phi_h\cdot\nabla w_l}{\lvert\nabla\phi_h\rvert},
\end{aligned}
\end{equation}
from which it follows that
\begin{equation}\label{e:dtfrac}
\der{}t1 \frac{\nabla(\phi_h + t w_l)}{\bigl\lvert\nabla(\phi_h + t w_l)\bigr\rvert}\Biggl|_{t=0}=
\frac{\lvert\nabla\phi_h\rvert\nabla w_l - \nabla\phi_h\frac{\nabla\phi_h\cdot\nabla w_l}{\lvert\nabla\phi_h\rvert}}{\lvert\nabla\phi_h\rvert^2}
= \left(1 - \frac{\nabla\phi_h}{\lvert\nabla\phi_h\rvert}\otimes\frac{\nabla\phi_h}{\lvert\nabla\phi_h\rvert}\right)\frac{\nabla w_l}{\lvert\nabla\phi_h\rvert}.
\end{equation}
From expressions~\eqref{e:nt} and~\eqref{e:dtfrac} and the fact that $\nabla\phi_h\big|_{\Gamma_{\mathrm p, h}}^{} = \partial_n\phi_h\big|_{\Gamma_{\mathrm p, h}}^{}$  almost everywhere on $\Gamma_{\mathrm p, h}$ follows that
\begin{equation}\label{e:n'}
\bs n' = \frac{P_T\nabla w_l}{\lvert\partial_n\phi_h\rvert} \qquad\text{a.\ e.\ on $\Gamma_{\mathrm p, h}$,}
\end{equation}
where
\begin{equation}
P_T = I - \bs n\otimes\bs n
\end{equation}
is the projector on the tangent plane of $\Gamma_{\mathrm p, h}$.
Substituting expression~\eqref{e:n'} into formula~\eqref{e:phi'}, we find that
\begin{equation}\label{e:phi'2}
\begin{aligned}
\psi'(q_h, p_h; \bs n) 
= 
\begin{cases}
\psi(q_h, p_h'; \bs n) &\text{on $\Gamma_{\mathrm w, h}\setminus\Gamma_{\mathrm p, h}$.}
\\[4pt] \displaystyle
\psi(q_h, p_h'; \bs n) 
- \alpha_V\Bigl(\frac{(P_T\nabla w_l)\cdot\nabla q_h}{\lvert\partial_n\phi_h\rvert}\pder{p_h}n1 + \pder{q_h}n1\frac{(P_T\nabla w_l)\cdot\nabla p_h}{\lvert\partial_n\phi_h\rvert}\Bigr)&\text{on $\Gamma_{\mathrm p, h}$,}
\end{cases}
\end{aligned}
\end{equation}
which, substituted in expression~\eqref{e:acfemdiff}, yields that $\forall q_h\in W_h$,
\begin{equation}\label{e:acfemdiff2}
\begin{aligned}
&\int_{\Omega_{\mathrlap h}} \nabla q_h\cdot\nabla p_h' - k^2\int_{\Omega_{\mathrlap h}} q_h p_h' 
- \int_{\Gamma_{\mathrlap{\mathrm p, h}}}\bigl(\nabla q_h\cdot\nabla p_h - k^2 q_h p_h\bigr)\frac{w_l}{\lvert\partial_{\!n}\,\phi_h\rvert}
+ \i k\int_{\Gamma_{\mathrlap{\text{out}}}} q_h p_h'
\\
&\qquad
+ \int_{\Gamma_{\mathrlap{\mathrm w, h}}} \psi(q_h,p_h'; \bs n)   
-\alpha_V\int_{\Gamma_{\mathrm p, h}}\Bigl(\frac{(P_T\nabla w_l)\cdot\nabla q_h}{\lvert\partial_n\phi_h\rvert}\pder{p_h}n1 + \pder{q_h}n1\frac{(P_T\nabla w_l)\cdot\nabla p_h}{\lvert\partial_n\phi_h\rvert}\Bigr)
\\
&\qquad
- \int_{\Gamma_{\mathrm p, h}}\! \pder{}n1\psi(q_h,p_h;\bs n) \frac{w_l}{\lvert\partial_{\!n}\,\phi_h\rvert}
- \sum_{S\in\mathscr S_h} \int\limits_{\Gamma_{\mathrm p, h}\cap S} \bs n^S\cdot\bigl\llbracket \psi(q_h, p_h;\bs n) \bs m\bigr\rrbracket\frac{w_l}{\lvert\partial_{n^{\!S}}\phi_h\rvert} 
+ \epsilon_h s_h(q_h, p_h')= 0,
\end{aligned}
\end{equation}
where we, for integral terms~3, 7 and~8, also have used that $w_l$ vanishes on $\partial\Omega_h\setminus\Gamma_{\mathrm p, h}$.

Now let $z_h\in W_h$ satisfy the adjoint equation
\begin{equation}\label{e:adjoint}
\begin{aligned}
\int_{\Omega_h}\!\!\nabla w_h\cdot\nabla z_h - k^2\int_{\Omega_h}\!\! w_h z_h + \i k\int_{\Gamma_{\mathrlap{\text{out}}}} w_h z_h
+ \int_{\Gamma_{\mathrm w, h}}\!\!\!\psi(w_h,z_h; \bs n) + \epsilon_h s_h(w_h,z_h) = \int_{\Gamma_{\mathrlap{\text{out}}}}w_h
\quad\forall w_h\in W_h.
\end{aligned}
\end{equation}
In particular, evaluating the adjoint equation for $w_h = p'_h$ yields that 
\begin{equation}\label{e:adjointph'}
\begin{aligned}
\int_{\Omega_h}\!\!\nabla p_h'\cdot\nabla z_h - k^2\int_{\Omega_h}\!\! p_h' z_h + \i k\int_{\Gamma_{\mathrlap{\text{out}}}} p_h' z_h
+ \int_{\Gamma_{\mathrm w, h}}\!\!\!\psi(p_h',z_h; \bs n) + \epsilon_h s_h(p_h',z_h) = \int_{\Gamma_{\mathrlap{\text{out}}}}p_h'.
\end{aligned}
\end{equation}
Substituting expression~\eqref{e:adjointph'} into expression~\eqref{e:acfemdiff2} with $q_h= z_h$,  the resulting expression reduces to 
\begin{equation}\label{e:acfemdiff3}
\begin{aligned}
&\int_{\Gamma_{\mathrlap{\text{out}}}}p_h'- \int_{\Gamma_{\mathrlap{\mathrm p, h}}}\bigl(\nabla z_h\cdot\nabla p_h - k^2 z_h p_h\bigr)\frac{w_l}{\lvert\partial_{\!n}\,\phi_h\rvert}
-\alpha_V\int_{\Gamma_{\mathrm p, h}}\!\Bigl(\frac{(P_T\nabla w_l)\cdot\nabla z_h}{\lvert\partial_n\phi_h\rvert}\pder{p_h}n1 + \pder{z_h}n1\frac{(P_T\nabla w_l)\cdot\nabla p_h}{\lvert\partial_n\phi_h\rvert}\Bigr)
\\
&\qquad 
- \int_{\Gamma_{\mathrm p, h}}\! \pder{}n1\psi(z_h,p_h;\bs n) \frac{w_l}{\lvert\partial_{\!n}\,\phi_h\rvert}
- \sum_{S\in\mathscr S_h} \int\limits_{\Gamma_{\mathrm p, h}\cap S} \bs n^S\cdot\bigl\llbracket \psi(z_h, p_h;\bs n) \bs m\bigr\rrbracket\frac{w_l}{\lvert\partial_{n^{\!S}}\phi_h\rvert}  = 0,
\end{aligned}
\end{equation} 
where we also have used that $s_h$ is symmetric and $\psi$ is symmetric in its first two arguments.

Finally, from expressions~\eqref{e:djp'} and~\eqref{e:acfemdiff3} we obtain the directional derivative expression
\begin{equation}
\begin{aligned}
dj(\phi_h; w_l ) &=  \int_{\Gamma_{\mathrlap{\mathrm p, h}}}\bigl(\nabla z_h\cdot\nabla p_h - k^2 z_h p_h\bigr)\frac{w_l}{\lvert\partial_{\!n}\,\phi_h\rvert}
\\
&\quad +\alpha_V\int_{\Gamma_{\mathrm p, h}}\!\Bigl(\frac{(P_T\nabla w_l)\cdot\nabla z_h}{\lvert\partial_n\phi_h\rvert}\pder{p_h}n1 + \pder{z_h}n1\frac{(P_T\nabla w_l)\cdot\nabla p_h}{\lvert\partial_n\phi_h\rvert}\Bigr)
\\
&\qquad 
+ \int_{\Gamma_{\mathrm p, h}}\! \pder{}n1\psi(z_h,p_h;\bs n) \frac{w_l}{\lvert\partial_{\!n}\,\phi_h\rvert}
+ \sum_{S\in\mathscr S_h} \int\limits_{\Gamma_{\mathrm p, h}\cap S} \bs n^S\cdot\bigl\llbracket \psi(z_h, p_h;\bs n) \bs m\bigr\rrbracket\frac{w_l}{\lvert\partial_{n^{\!S}}\phi_h\rvert},
\end{aligned}
\end{equation}
under perturbation~\eqref{e:phiht}, where $p_h$ solves state equation~\eqref{e:acfemperturb} for $t=0$ and $z_h$ adjoint equation~\eqref{e:adjoint}.
\end{appendices}

\clearpage
\bibliographystyle{plain}
\bibliography{3DPhasePlugOpt}

\end{document}